\documentclass[onefignum,onetabnum]{siamart190516}
\usepackage{amsmath,amsfonts,amssymb,mathrsfs,cases,mathdots,colordvi,url}
\usepackage{enumerate}
\usepackage{amssymb,latexsym}
\usepackage{graphicx}
\usepackage{amstext}
\usepackage{amsopn}
\usepackage{amsbsy}
\usepackage{amscd}
\usepackage{amsxtra}
\usepackage{graphics}
\usepackage{booktabs,multirow}
\usepackage{tabularx}
\usepackage{makecell}
\usepackage{bbding}
\usepackage{pifont}
\usepackage{hyperref}
\usepackage{tikz}   
\usetikzlibrary{shapes.geometric,arrows} 
\tikzstyle{block} = [rectangle, draw, fill=white!50,   
text width=33em, text centered, rounded corners, minimum height=3em] 
\tikzstyle{line} = [draw, -latex']  
\tikzstyle{arrow} = [thick,->,>=stealth]

\newtheorem{example}{Example}[section]

\usepackage{lipsum}
\usepackage{amsfonts}
\usepackage{graphicx}
\usepackage{epstopdf}
\usepackage{algorithmic}
\ifpdf
  \DeclareGraphicsExtensions{.eps,.pdf,.png,.jpg}
\else
  \DeclareGraphicsExtensions{.eps}
\fi

\newsiamremark{remark}{Remark}
\newsiamremark{hypothesis}{Hypothesis}
\crefname{hypothesis}{Hypothesis}{Hypotheses}
\newsiamthm{claim}{Claim}

\headers{Generic Property of Partial Calmness}{R. Ke, W. Yao, J.J. Ye, and J. Zhang}

\title{Generic property of the partial calmness condition for bilevel programming problems\thanks{Submitted to the editors DATE.
\funding{Ke's work is supported by National Science Foundation of China 71672177. The research of Ye is supported by NSERC. Zhang's work is supported by the Stable Support Plan Program of Shenzhen Natural Science Fund (No. 20200925152128002), National Science Foundation of China 11971220, Shenzhen Science and Technology Program (No. RCYX20200714114700072). The alphabetical order of the authors indicates the equal contribution to the paper.}}}

\author{Rongzhu Ke\thanks{School of Economics and Research Center for China's Digital Economy, Zhejiang University, Hangzhou,
		Zhejiang 310058, China (\email{rzke455@zju.edu.cn}). } 
	\and Wei Yao\thanks{Department of Mathematics, Southern University of Science and Technology, Shenzhen, Guangdong, China (\email{yaow@sustech.edu.cn}). }
	\and Jane J. Ye\thanks{Department of Mathematics
		and Statistics, University of Victoria, Victoria, B.C., Canada V8W 2Y2 (\email{janeye@uvic.ca}).}
	\and Jin Zhang\thanks{Corresponding author. Department of Mathematics, Southern University of Science and Technology, and National Center for Applied Mathematics Shenzhen, Guangdong, China
		(\email{zhangj9@sustech.edu.cn}).} 
}

\usepackage{amsopn}

\ifpdf
\hypersetup{
  pdftitle={Generic Property of Partial Calmness},
  pdfauthor={R. Ke, W. Yao, J.J. Ye, and J. Zhang}
}
\fi

\begin{document}

\maketitle

\begin{abstract}
  The partial calmness for the bilevel programming problem (BLPP) is an important condition which ensures that a local optimal solution of BLPP is a local optimal solution of a partially penalized problem where the lower level optimality constraint is moved to the objective function and hence a weaker constraint qualification can be applied. In this paper we propose a sufficient condition in the form of a partial error bound condition which guarantees the partial calmness condition. We analyse the partial calmness for the combined program  based on the Bouligand (B-) and the Fritz John (FJ) stationary conditions from a generic point of view. Our main result states that the partial error bound condition  for the combined programs based on B and FJ conditions are generic for an important setting with applications in economics and hence the partial calmness for the combined program is not a particularly stringent assumption. Moreover we derive optimality conditions for the combined program for the generic case without any extra constraint qualifications and show the exact equivalence between our optimality condition and the one by Jongen and Shikhman given in implicit form. Our arguments are based on Jongen, Jonker and Twilt's generic (five type) classification of the so-called generalized critical points for one-dimensional parametric optimization problems and Jongen and Shikhman's generic local reductions of BLPPs.
\end{abstract}



\section{Introduction}
\label{intro}

In this paper we consider the following bilevel programming problem (BLPP):
\begin{equation}\label{BLPP}\tag{BLPP}
\min_{x,y}\  F(x,y)\quad
\mathrm{s.t.}\quad y\in S(x), \  G(x,y)\leq0, 
\end{equation}
where $S(x)$ denotes the solution set of the lower level program
\begin{equation}\label{Lx}\tag*{$L(x)$}
\min_{y} \,f(x,y) \quad 
\mathrm{s.t.}
\quad  g(x,y)\leq0. 
\end{equation}
For convenience, we denote the feasible set of the lower level program by
$
Y(x):=\big\{y\in\mathbb{R}^m: g(x,y)\leq0\big\}.
$
Here $x\in\mathbb{R}^n$ and the mappings  $F,f:\mathbb{R}^n\times\mathbb{R}^m\rightarrow\mathbb{R}$, 
$G:\mathbb{R}^n\times\mathbb{R}^m\rightarrow\mathbb{R}^q$, 
$g:\mathbb{R}^n\times\mathbb{R}^m\rightarrow\mathbb{R}^p$. We assume that $F, G, f, g$ are continuously differentiable unless otherwise specified.

BLPPs have broad applications. For example, they appear in the principal-agent moral hazard problem \cite{M7599}, electricity markets and networks \cite{BJ05}. Recently, BLPP has been applied to hyperparameters optimization and and meta-learning in machine learning; see e.g. \cite{KBHP08,FFSGP18,LMYZZ20,YYZZ}. More applications can be found in the monographs \cite{SIB97,B98,D18,DKPK15}, the survey on bilevel optimization \cite{Dempebook} and the references therein. 

The main difficulty in studying BLPPs is that the constraint $y\in S(x)$ is a global one. The classical and most often used approach to study BLPPs is to replace the lower level problem by the Karush-Kuhn-Tucker (KKT) condition and minimize over the original variables as well as multipliers. The resulting problem is the so-called mathematical program with complementarity/equilibrium constraints (MPCC/MPEC) (\cite{MPEC1,MPEC2}), which itself is a difficult problem if treated as a nonlinear programming problem (NLP) since the Mangasarian Fromovitz constraint qualification (MFCQ) is violated at any feasible point of MPEC, see Ye et al. \cite{YZZ97}. In general, this approach is only applicable to BLPPs where the lower level problem is a convex optimization problem in variable $y$ since KKT condition is not sufficient for $y\in S(x)$ when the lower level program is not convex.  As a matter of fact, it was pointed out by Mirrlees \cite{M7599} that an optimal solution of BLPP may not even be a stationary point of the resulting single-level problem by  KKT/MPEC approach (usually called the first order approach in economics).  Moreover, since the reformulated problem   involves additionally the  multipliers of the lower level problem as variables, as pointed out by Dempe and Dutta in \cite{DD12}, it is possible that $(\bar{x},\bar{y},\bar{u})$ is a local optimal solution of the reformulated problem but $(\bar{x},\bar{y})$ is not a local optimal solution of BLPP even when the lower level problem is convex. 

To deal with BLPPs without convexity assumption on the lower level problem, the value function approach was proposed by Outrata \cite{O90} for numerical purpose and used by Ye and Zhu \cite{YZ95} for optimality conditions. By defining the value (or the marginal) function as
$
V(x):=\inf_y\big\{f(x,y): g(x,y)\leq0\big\},
$
one can replace the original BLPP by the following
equivalent problem:
\begin{equation}\label{VP}\tag{VP}
\min_{x,y}\  F(x,y)\quad
\mathrm{s.t.}\quad
f(x,y)-V(x)\leq0,\ 
g(x,y)\leq0,\ 
G(x,y)\leq0.
\end{equation}
The resulting problem \eqref{VP}  is clearly equivalent to the original bilevel problem. However, since the value function constraint is actually an equality constraint, the nonsmooth MFCQ for \eqref{VP} will never hold (cf. \cite[Prop. 3.2]{YZ95}, \cite[Section 4]{Ye20}). To derive necessary optimality conditions for BLPPs, Ye and Zhu \cite{YZ95} proposed the partial calmness condition for \eqref{VP} under which the difficult constraint $f(x,y)-V(x)\leq0$ is penalized to the objective function. Under the partial calmness condition for \eqref{VP}, the usual constraint qualifications such as MFCQ can be satisfied for the partially penalized (VP). Consequently, by using a nonsmooth multiplier rule, one can derive a KKT type optimality condition for \eqref{VP} provided that the value function is Lipschitz continuous. We refer the reader to \cite{YZ95,Ye04,Ye06,DDM07,DZ13,MNP12} for more discussions.

Although it was proved in \cite{YZ95} that the partial calmness condition for \eqref{VP} holds automatically for the minmax problem and the bilevel program where the lower level problem is linear in both variables $x$ and $y$, the partial calmness condition for \eqref{VP} has  been shown to be a celebrated but restrictive assumption (see \cite{HS11,Ye20,MMZ20}). To improve the value function approach, Ye and Zhu \cite{YZ10} proposed   a combination of the classical KKT and the value function approach. The resulting problem is  the combined problem using  KKT condition: 
\begin{equation}\label{CP}\tag{CP}
\begin{aligned}
\min_{x,y,u}\ & F(x,y) \\
\mathrm{s.t.\ \,} &f(x,y)-V(x)\leq0, \ 
\nabla_y f(x,y)+u\nabla_y g(x,y)=0, \\
&g(x,y)\leq0,\ u\geq0,\ u^T g(x,y)=0,
\ G(x,y)\leq0, 
\end{aligned}
\end{equation}
where $u\nabla_y g(x,y):=\sum_{i=1}^p u_i\nabla g_i(x,y)$. Similarly to \cite{YZ95},  to deal with the fact that  the nonsmooth MFCQ also fails for \eqref{CP}, the following partial calmness condition for \eqref{CP} was proposed in  \cite{YZ10}.
\begin{definition}[\bf Partial calmness for \eqref{CP}]
	Let $(\bar{x},\bar{y},\bar u)$ be a feasible solution of \eqref{CP}. We say that \eqref{CP} is partially calm at $(\bar{x},\bar{y},\bar{u})$ if there exists $\mu\geq 0$ such that $(\bar{x},\bar{y},\bar{u})$ is a local solution of the following partially penalized problem:
	\begin{equation}\label{CP1}\tag{\text{CP$_\mu$}}
	\begin{aligned}
	\min_{x,y,u}\ & F(x,y)+\mu\big(f(x,y)-V(x)\big) \\
	\mathrm{s.t.\ \,}
	& \nabla_y f(x,y)+u\nabla_y g(x,y)=0, \\
	& g(x,y)\leq0,\ u\geq0,\ u^T g(x,y)=0,
	\ G(x,y)\leq0. 
	\end{aligned}	
	\end{equation} 	
\end{definition}
Since there are more constraints in \eqref{CP} than in \eqref{VP}, the partial calmness for \eqref{CP} is a weaker condition than the one for \eqref{VP}. Moreover, it was shown in \cite{YZ10,Ye11} that under the partial calmness condition and various MPEC variant constraint qualifications for \eqref{CP1}, a local optimal solution of \eqref{CP} must be a corresponding MPEC variant stationary point based on the value function provided the value function is Lipschitz continuous. See \cite{YZ10,Ye11,Ye20} for more details. In the same spirit, Nie et al. \cite{nie2017bilevel}  proposed a numerical algorithm for globally solving polynomial bilevel programs by adding the generalized critical point condition for the lower level program. Recently, Nie et al. \cite{nie2020bilevel} expressed the Lagrange multipliers of the lower level polynomial problem as a function of the upper and lower level variables and hence eliminated the multiplier variables in the reformulation of BLPPs in \cite{nie2017bilevel}.

The reformulation \eqref{CP}, however, requires the existence of KKT condition at each optimal solution of the lower level program (see \cite[Proposition 3.1]{YZ10}). In order to deal with the case where KKT condition may not hold at all the solutions of the lower level program, we propose the following combined program using the  Fritz John (FJ) condition\footnote{ FJ condition was firstly used by Allende and Still \cite{AS13} for BLPPs, where they replaced the lower level program by FJ condition and the resulting problem becomes \eqref{CPFJ} without the value function constraint.}:
\begin{equation}\label{CPFJ}\tag{CPFJ}
\begin{aligned}
\min_{x,y,u_0, u}\ & F(x,y) \\
\mathrm{s.t.\ \,} &f(x,y)-V(x)\leq0, \\
&u_0\nabla_y f(x,y)+u\nabla_y g(x,y)=0,\ g(x,y)\leq0,\  (u_0,u) \geq0, \\
&u^T g(x,y)=0, \ {\|(u_0,u)\|_1=1},
\ G(x,y)\leq 0,
\end{aligned}
\end{equation}
where {the one-norm $\|(u_0,u)\|_1:=\sum_{i=0}^p |u_i|=\sum_{i=0}^p u_i$ when $(u_0,u) \geq0$.}
A relationship between problems  \eqref{CPFJ} and \eqref{BLPP} can be established by similar arguments as in \cite[Proposition 3.1]{YZ10}. 
Similarly to \cite{YZ10}, we propose the partial calmness condition for \eqref{CPFJ}. 
\begin{definition}[\bf Partial calmness for \eqref{CPFJ}]\label{Defn1.2}
	Let $(\bar{x},\bar{y},\bar u_0,\bar{u})$ be a feasible solution of \eqref{CPFJ}. We say that \eqref{CPFJ} is partially calm at $(\bar{x},\bar{y},\bar u_0,\bar{u})$ if there exists $\mu\geq 0$ such that $(\bar{x},\bar{y},\bar u_0,\bar{u})$ is a local solution of the partially penalized problem:
	\begin{equation}\label{CP2}\tag{\text{CPFJ$_\mu$}}
	\begin{aligned}
	\min_{x,y,u_0, u}\ & F(x,y)+\mu\big(f(x,y)-V(x)\big) \\
	\mathrm{s.t.\ \,}
	&u_0\nabla_y f(x,y)+u\nabla_y g(x,y)=0,\ g(x,y)\leq0,\ (u_0,u)\geq0, \\
	&u^T g(x,y)=0,\ {\|(u_0,u)\|_1=1},
	\ G(x,y)\leq0. 
	\end{aligned}
	\end{equation} 	
\end{definition}

Inspired by \cite{AHO18,GY17,GY19,Ye20} where the authors study the mathematical program with a generalized equation involving the regular normal cone, we consider the combined program with the Bouligand (B-) stationary condition for the lower level program:
\begin{equation}\label{CPB}\tag{CPB}
\begin{aligned}
\min_{x,y}\ & F(x,y)\\
\mathrm{s.t.\,}&f(x,y)-V(x)\leq0, \\
&0\in \nabla_y f(x,y)+\widehat{N}_{Y(x)}(y),
\ G(x,y)\leq0.
\end{aligned}
\end{equation}
Here $\widehat{N}_{Y(x)}(y)$ is the regular normal cone to $Y(x)$ at $y$ if $y\in Y(x)$ (see Definition \ref{TNcone}) and is equal to the empty set if $y\not \in Y(x)$.  
The combined problem \eqref{CPB} is equivalent to the original bilevel program since the condition  
\begin{equation}\label{ge}
0\in \nabla_y f(x,y)+\widehat{N}_{Y(x)}(y)
\end{equation} 
is a necessary optimality condition for $y\in S(x)$.  Indeed, if $y$ is a solution to the lower level problem $L(x)$, then  there is no descent direction starting from $y$ and hence B-condition holds which means that \eqref{ge}  must hold; see Section~\ref{sec2} for more details. Similarly, we propose the following partial calmness condition for \eqref{CPB}. 
\begin{definition}[\bf Partial calmness for \eqref{CPB}]\label{ourcalmness}
	Let $(\bar x, \bar{y} )$ be a feasible solution of \eqref{CPB}.
	We say that \eqref{CPB} is partially calm at $(\bar{x},\bar{y})$ if there exists $\mu\geq 0$ such that $(\bar{x},\bar{y})$ is a local solution of the following partially penalized problem: 
	\begin{equation}\label{CPB1}\tag{\text{CPB$_\mu$}}
	\begin{aligned}
	\min_{x,y}\ & F(x,y)+\mu\big(f(x,y)-V(x)\big) \\
	\mathrm{s.t.\ }
	&0\in \nabla_y f(x,y)+\widehat{N}_{Y(x)}(y),\ G(x,y)\leq0.
	\end{aligned}
	\end{equation}	
\end{definition} 
Note that there are other alternative combined programs. For example, we may consider the combined program
\begin{equation}\label{CP4}\tag{\text{CP$_{FJ}$}}
\begin{aligned}
\min_{x,y}\ & F(x,y) \\
\mathrm{s.t.\ \,} &f(x,y)-V(x)\leq0, \\
&(x,y) \in \Sigma_{FJ},\ G(x,y)\leq 0,
\end{aligned}
\end{equation}
where $\Sigma_{FJ}$ is the set of FJ points defined as in the Fritz John condition \eqref{FJ} and the partial calmness condition can be similarly defined.

From the discussion above, we can define various types of combined programs and the corresponding partial calmness conditions based on different stationary conditions for the lower level program. But combined programs can be classified into two types. One  involves  multipliers as extra variables such as \eqref{CP} and \eqref{CPFJ} and the other does not involve any extra variables such as \eqref{CPB} and  \eqref{CP4}.  These two types of combined programs have its own advantages and disadvantages. The one with multipliers such as \eqref{CP} and \eqref{CPFJ} are explicit and easier to handle while the one without multipliers are implicit and harder to handle. An advantage of  using a combined program without  multipliers such as \eqref{CPB} is that the constraint qualification for \eqref{CPB1}  is in general weaker than the one for \eqref{CP1}. By {a counterexample in \cite[Appendix]{AHO18}}, it is possible that  a feasible solution of  \eqref{CPB1} satisfies the calmness condition at $(\bar x,\bar y, 0)$ while the one for the corresponding \eqref{CP1} does not satisfy the calmness condition at $(\bar x,\bar y, u, 0)$ for any multiplier $u$. By \cite[Proposition 5]{GY17}, we see that for the case where $Y(x)$ is independent of $x$, if the lower level multipliers are not unique, then the MPEC-MFCQ never holds for \eqref{CP1} and  even the very weak constraint qualification MPCC-GCQ can be violated for \eqref{CP1} while the calmness condition can hold at $(\bar x,\bar y, 0)$ for \eqref{CPB1}.

\textbf{Questions.} The main goal of this paper is to investigate two questions: 
\begin{equation}\label{q1}\tag{Q1}
\text{How stringent is the partial calmness condition for combined program? }
\end{equation}
This question is at present far from being solved. It is not clear, a priori, whether the partial calmness for the combined program is a typical property for BLPPs. 

The second question is: 
\begin{equation}\label{q2}\tag{Q2}
\begin{aligned}
&\text{Is there a generic sufficient condition to ensure}\\
&\quad\text{ the partial calmness condition for the combined program? }
\end{aligned}
\end{equation}
Here a property is said to be generic if it holds for all problems generated by data in an open and dense subset or, more generally, in the intersection of countably many open and dense subsets. It is well known that the nondegeneracy of KKT points is a generic property for standard nonlinear programming problems (cf. \cite{JJT00}) and it has been also shown that a MPCC variant  of Linearly Independence Constraint Qualification (LICQ) called MPEC-LICQ is a generic property for MPCC/MPEC (cf. \cite{SS01,AS13}). {The partial calmness condition plays an important role in the theory and applications of BLPPs, particularly in deriving optimality conditions since the partially penalized problem does not have the difficult value function constraint in it;  \cite[Chapter 6]{M18}, \cite{M20}, and \cite{Ye20}. }

\textbf{Contributions.} To answer the above questions, by virtue of  the exact penalty principle of Clarke, we propose a sufficient condition for the partial calmness in the form of a partial error bound condition in Definition \ref{vfcq}. In the case where the upper level variable has dimension one, our main result (Theorems \ref{thm1}) shows that the partial error bound (PEB) and consequently the partial calmness condition for \eqref{CP4} and \eqref{CPB} are generic for BLPPs. Moreover, by an equivalence result in Theorem \ref{thmequiv}, we have that the partial calmness condition for \eqref{CPFJ} is also generic. This   gives an affirmative answer to \eqref{q2} in this case.  Hence, generically, all local optimal solutions of BLPP are local optimal solutions of \eqref{CPB1} for some nonnegative constant $\mu\geq 0$. Similarly, generically, all local optimal solutions to \eqref{CPFJ} are local optimal solutions of \eqref{CP2} for some nonnegative constant $\mu\geq 0$. Combining these results and recent results in \cite{GY17,GY19} for mathematical program with a parametric generalized equation gives some  new sharp necessary optimality conditions for BLPPs. For the generic cases,  we are able to show that MPEC-LICQ holds for partially perturbed problem  \eqref{CP2} without the upper level constraint. Consequently in Theorem \ref{condition1}, we derive necessary optimality conditions for the generic cases explicitly in problem data and without any extra constraint qualifications. Further, in Theorem~\ref{eqivthm}, we prove that our optimality conditions for the generic case are equivalent to the one derived by Jongen and Shikhman in \cite{JS12} and hence partially answer a question raised in \cite[Remark 4.2]{JS12}. When there is no upper and lower level constraints in BLPPs, our result (Corollary~\ref{thmno}) shows that the condition in  \cite[Theorem 2.1]{YZ10} is a generic one. This gives a rigorous proof of Mirrlees' arguments about $f$-critical points in \cite{M7599}. 

\textbf{Technique.} By adopting a parametric programming point of view, our arguments to prove the main result (Theorems \ref{thm1}) are based on Jongen, Jonker and Twilt's generic (five type) classification of the so-called generalized critical points (g.c.) for one-dimensional parametric optimization problems \cite{JJT86}, and Jongen and Shikhman's generic local reductions of BLPPs in the case where the dimension of the upper variable is one \cite{JS12}. To obtain local reduction, Jongen and Shikhman \cite{JS12} studied the structure of points being local minimizers for \ref{Lx}. In this paper, to study the combined problems \eqref{CP}, \eqref{CPFJ}, \eqref{CPB} and \eqref{CP4},  we study the structure of points being KKT, FJ, and B-points for \ref{Lx}, by using a similar analysis in \cite{JS12}. The roadmap in Figure \ref{roadmap} is helpful at first glance and the details are left to  Sections~\ref{sec5},\ref{sec6}, \ref{sec7}. Note that although the classification is done under the assumption that the parameter is one {dimensional}, it gives some  indications on the situation in the general case where the parameter is from a multi-dimensional space; see \cite{JJT00,DGJ09,JS12,DJS} for more discussions.

The generic classification of g.c. points is studied within the scope of singularity theory (see e.g. \cite{BL75}). The proofs of results in \cite{JJT86} are mainly based on the   transversality theory. We refer the reader to \cite{JJT00} for a detailed account of the topological aspects of nonlinear optimization. It has proved that topology method is a powerful tool for studying the fine structure of optimization problems (e.g. \cite{JJT00,JS12,SS01,AS13,DJS} and the references given there). This work is intended as an attempt to further motivate the study of topological aspects of BLPPs. And it is of interest to bring together two areas in studying BLPPs. One is to study the generic structure of bilevel problems in a neighborhood of its solution set, cf. \cite{DGJ09,JS12,DJS}. Another one is to study constraint qualifications and optimality conditions for BLPPs, see e.g. \cite{YZ10,M20,Ye20}. 

\textbf{Outline.} The remaining part of the paper is organized as following. In Section \ref{sec2}, we gather some preliminaries in variational analysis. Some motivational applications in economics will be described in Section \ref{application}. In Section \ref{sec5}, we introduce the main technique and state our main result whose proof is given in Section \ref{sec6}. As a consequence, we study optimality conditions for BLPPs in Section \ref{sec7}.

\textbf{Symbols and Notations.} Our notation is basically standard. By $0_m$, we mean it is the zero vector in $\mathbb{R}^m$. For a vector $z\in \mathbb{R}^d$, we denote its Euclidean norm by $\|z\|$. {Simply by $\mathbb{B}$, we denote the closed unit ball centered at the origin of the space in question.}
For a matrix $A$, we denote by $A^T$ its transpose. The inner product of two vectors $x,y$ is denoted by $x^T y$ or $\langle x,y\rangle$. For $z\in\mathbb{R}^d$ and  $\Omega\subseteq\mathbb{R}^d$, we denote by $\mathrm{dist}_\Omega(z)$ the distance from $z$ to $\Omega$. By $z\stackrel{\Omega}{\rightarrow}\bar{z}$ we mean $z\rightarrow \bar z$ and $z\in \Omega$. The polar cone of a set $\Omega$ is $\Omega^\circ:=\{v: \langle v, z\rangle\leq0\ \forall z\in\Omega\}$. We also denote by $\mathrm{co}\,\Omega$ the convex hull of $\Omega$. For a function $h:\mathbb{R}^d\rightarrow\mathbb{R}$, we denote the gradient vector and the Hessian of $h$ at $z$ by $\nabla h$ and $\nabla^2 h$, respectively. {Given a differentiable one-dimensional function $\Psi$, $D\Psi$ denotes its derivative.} Let $\Phi$ be a set-valued map. We denote its graph by $\mathrm{gph}\Phi:=\{(z,w): w\in\Phi(z)\}$. Let $C^3(\mathbb{R}^d)$ denotes the space of three times continuously differentiable functions on $\mathbb{R}^d$. We denote by $C_s^3$ the strong (or Whitney) $C^3$-topology and refer the reader to \cite[Section 6.1]{JJT00} for more details on $C^3_s$-topology. 

\section{Preliminaries on variational analysis}\label{sec2}

In this section, we gather some preliminaries in variational analysis that will be needed. The reader is referred to Rockafellar-Wets \cite{RW98} and the references therein for more discussions on the subject. 

First we recall the definition for the tangent cone and the normal cones. 

\begin{definition}[\bf Tangent  and normal cones]\label{TNcone}
	\cite[Definitions 6.1 and 6.3]{RW98}
	Given a set $\Omega\subseteq\mathbb{R}^d$ with $\bar{z}\in\Omega$, the (Bouligand-Severi) tangent/contingent cone to $\Omega$ 
	at $\bar{z}$ is a closed cone defined by
	\begin{equation*}
	T_\Omega(\bar{z}):=\limsup_{t\downarrow0}\frac{\Omega-\bar{z}}{t}=
	\Big\{u\in\mathbb{R}^d: \exists\,t_k\downarrow 0, u_k\rightarrow u\ \mathrm{with}\ \bar{z}+t_k u_k\in\Omega,\ \forall\,k\Big\}.
	\end{equation*}
	The (Fr\'{e}chet) regular normal cone and the (Mordukhovich) limiting/basic normal
	cone to $\Omega$ at $\bar{z}\in\Omega$ are closed cones defined by
	$
	\widehat{N}_\Omega(\bar{z}):=\Big\{v\in\mathbb{R}^d: \limsup_{z\stackrel{\Omega}{\rightarrow}\bar{z}}\frac{\langle v,z-\bar{z}\rangle}{\|z-\bar{z}\|}\leq0\Big\}
	$
and 
$
	\quad N_\Omega(\bar{z}):=\left\{v\in\mathbb{R}^d: \exists\,z_k\stackrel{\Omega}{\rightarrow}\bar{z},\ v_k\in \widehat{N}_\Omega(z_k)\ \mathrm{with}\ v_k\rightarrow v\ \mathrm{as}\ k\rightarrow\infty\right\}.
	$
\end{definition}
When the set $\Omega$ is convex, the tangent/contingent cone and the regular/limiting normal cone reduce to the derivable  cone and the normal cone in the sense of convex analysis, respectively; see e.g., \cite[Theorem 6.9]{RW98}.

It is well-known that a necessary condition for $\bar z$ to be locally optimal to the problem of minimizing a differentiable function $f_0$ over a set $\Omega \subseteq \mathbb{R}^d$ is 
\begin{equation} 
\langle \nabla f_0(\bar z), u\rangle \geq 0 \quad \mbox{ for all } u\in T_\Omega(\bar z).
\label{decentd} 
\end{equation}
By the tangent-normal polarity (\cite[Theorem 6.28]{RW98}) 
$\widehat{N}_\Omega(\bar{z})=T_\Omega(\bar{z})^\circ$, the condition (\ref{decentd}) is the same as 
\begin{equation}  -\nabla f_0(\bar z) \in \widehat{N}_\Omega(\bar{z}) \ \ \mbox{ or } \ \ 0\in \nabla f_0(\bar z)+ \widehat{N}_\Omega(\bar{z}) . \label{Bcondition}
\end{equation}
Since the condition (\ref{decentd}) is called the Bouligand (B-) stationary condition in the literature, and conditions (\ref{decentd}) and (\ref{Bcondition})  are equivalent, we also call the condition  (\ref{Bcondition}) B-condition and any $y$ satisfying B-condition a B-point in this paper. 

We now review some concepts of stability of a set-valued map. 

\begin{definition}[\bf Calmness \cite{YY97,RW98}]
	Let $\Gamma:\mathbb{R}^d\rightrightarrows\mathbb{R}^s$ be a set-valued map and $(\bar{z},\bar w)\in \mathrm{gph}\Gamma$. We say that $\Gamma$ is calm (or pseudo upper-Lipschitz continuous) at $(\bar{z},\bar w)$ or equivalently its inverse map $\Gamma^{-1}$ is metrically subregular at $(\bar w,\bar{z})$ if there exist a neighbourhood 
	$\mathcal{Z}$ of $\bar{z}$, a neighborhood $\mathcal{W}$ of $\bar w$ and $\kappa\geq0$ such that		
	\begin{equation*}
	\Gamma(z)\cap \mathcal{W}\subseteq\Gamma(\bar{z})+\kappa\|z-\bar{z}\|\mathbb{B},\ \forall z\in\mathcal{Z}.
	\end{equation*}		
\end{definition}

\begin{definition}[\bf Metric subregularity constraint qualification] 
	
	Let $g(\bar x,\bar y)\leq 0$, where $g:\mathbb{R}^n\times\mathbb{R}^m\rightarrow\mathbb{R}^p$. We say that the system $g(x,y)\leq 0$ satisfies the metric subregularity constraint qualification (MSCQ) at $(\bar x,\bar y)$ if the perturbed set-valued map 
	$\Gamma( z): =\{ (x,y): g(x,y)\leq z \}$ is calm at $(0,\bar x,\bar y)$ or equivalently the set-valued map $\Gamma^{-1}(x,y):=-g(x,y)+\mathbb{R}^p_-$ is metrically subregular at $(\bar x,\bar y, 0)$.
\end{definition}

It is known that the linear CQ (when $g$ is affine) and  MFCQ both imply MSCQ.  A summary of more CQs which implies MSCQ is given in \cite[Theorem 7.4]{YZ18}.

\section{Motivating applications in economics}\label{application}

Our research is motivated by a number of important applications in economics. Below we describe some of them. 

\textbf{$\bullet$ The classic principal-agent moral hazard problem}

In the situation of {the moral hazard principal-agent problem}, the agent chooses an action $a\in A$ that is unobservable to a principal. The action $a$ influences the random outcome $s\in S$ through a parameterized probability measure $P(\cdot|a)$ (see \cite{M7599}). The outcome $s$ leads to a monetary payoff  $\pi(s)$, which accrues directly to the principal. Since the principal cannot monitor the agent's action but only the outcome, the principal will pay the agent conditional on the observed outcome. Thus the compensation contract specifies a wage of $w(s)$. If $S$ is a discrete subset of $\mathbb{R}$ with cardinality $k$, then the contract is a vector $w=(w_1,\dots,w_k)\in\mathbb{R}^k$. The principal's utility for outcome $s$ is a function of the net value $\pi(s)-w(s)$ and is denoted by $v(\pi(s)-w(s))$. The agent receives the wage $w(s)$ and pays the action $a$, then he receives utility $u(w(s),a)$. Hence their expected utility functions are 
\begin{equation*}
V(w,a)=\int_S v(\pi(s)-w(s)) dP(s|a),\quad
\mathrm{and}\quad U(w,a)=\int_S u(w(s),a) dP(s|a).
\end{equation*}
Therefore, the principal's problem in parameterized distribution terms can be stated:  
\begin{align}\label{pa}
&\max_{w,a} V(w,a)\tag{PA}\\
\mathrm{s.t.}\ &a\in\mathrm{arg}\max_{a'\in A} U(w,a'),\label{GIC}\tag{IC}\\
&U(w,a)\geq\underline{U}.\label{GIR}\tag{IR}
\end{align}
The first constraint is the incentive compatibility constraint which means that the agent will only take actions that maximize his own expected utility. Here we consider optimistic formulation. That is, we assume implicitly that the agent will choose the action most beneficial to the principal if he is indifferent between several different actions. The second constraint is the individual rationality constraint (or participation constraint) which guarantees that the agent earns at least his reservation utility $\underline{U}$.

The difficulty to solve the above moral hazard problem is well-recognized (e.g., Conlon \cite{C09}). But according to our method developed in this paper, two important categories of problem can fall into the special class of genericty. The first special case is binary outcome distribution $k=2$. The second special case is that the contract takes a linear form $w(s)=\alpha +\beta s$. 

\textbf{$\diamond$ Binary Outcome}

Binary outcome model is useful in economic literature (cf. Grossman
and Hart \cite{GH}).\footnote{%
	For example, in the case of social insurance, the outcome of job search is a
	binary outcome, unemployment or employment. In the case of coroperate
	finance, the outcome of investment project could be a binary outcome,
	success or not success. } 
Suppose that the principal is risk-neutral (i.e., $v(r)=r$) and assume that agent's utility is additively separable (i.e., {$u(w,a)=u(w)-c(a)$}). Consider the case $S=\{s_{1},s_{2}\}$ and assume that $u$ is increasing. It is known that at  the optimal solution, the IR constraint must be binding {\cite[Proposition 2]{GH}}.  Hence by setting 
\begin{equation*}
p_{i}(a):=P(\{s_{i}\}|a),\ \pi _{i}:=\pi (s_{i}),\ w_{i}:=u(w(s_{i})),\ \mathrm{and}\ x:=w_{1}-w_{2},\ y:=a,
\end{equation*}%
the principal-agent problem becomes 
\begin{equation*}
\min_{x,y}F(x,y) 
\quad \mathrm{s.t.}\quad y\in \mathrm{arg}\min_{y^{\prime }\in A} \big\{-p_{1}(y^{\prime
})x+c(y^{\prime })\big\},
\end{equation*}%
where 
\begin{equation*}
F(x,y):=p_{1}(y)\left[ h\big(\underline{U}+c(y)+p_{2}(y)x\big)-\pi _{1}
\right] 
+p_{2}(y)\left[ h\big(\underline{U}+c(y)-p_{1}(y)x\big)-\pi _{2}
\right],
\end{equation*}
and $h(\cdot ):=u^{-1}$($\cdot $) is increasing. This shows that the case $k=2$ in \eqref{pa} becomes the case $n=1$ in \eqref{BLPP}. 

\textbf{$\diamond$ Linear Contract}

Linear contracts play an important role in contract theory. Not only the
linear contract is intuitive and simple in practice, but also Holmstr\"{o}m
and Milgrom \cite{HM87} provided a theoretical foundation for the optimality
of linear contract in a dynamic setting where the agent with constant
absolute risk aversion controls the drift of a Brownian motion. 

Assume that the agent's utility is additively separable, we consider a
linear contract $w(s)=\alpha +\beta s$, where the slope $\beta $ is the
\textquotedblleft incentive intensity". Then agent's expected utility
function becomes 
$
\int_{S}u(\alpha +\beta s)dP(s|a)-c(a).
$
Assuming that  $u$ is an increasing function, we  can solve $\alpha =\alpha (\beta ,a)$ in terms of $\beta $ and $a$ by the
binding IR constraint. Therefore, relabel
the notation 
$x=\beta, y=a,$
for every implementable $y=a$ (i.e., there exists some $x$ such that $(x,y)$
is bilevel feasible), we can write the original moral hazard problem as 
\begin{align*}
& \min_{x}-\int_{S}v\Big(\pi (s)-\big(\alpha (x,y)+xs\big)\Big)dP(s|y) \\
\mathrm{s.t.\ }& y\in \mathrm{arg}\min_{y^{\prime }\in A}\Big\{c(y^{\prime
})-\int_{S}u(\alpha (x,y)+xs)\,dP(s|y^{\prime })\Big\},
\end{align*}%
which belongs to the class specifed by \eqref{BLPP}. \footnote{%
	Note that here we do not include the lower level variable $y$ into the upper
	level optimization problem because $y$ is given. So we obtain a
	characterization of $x$ for every implementable $y$, which is slightly
	different. However, if the agent is risk-neutral or the contract is chosen
	to make $u(w(s))=\alpha +\beta s$, we can obtain an exactly form of problem
	specified by \eqref{BLPP}.}

\textbf{$\bullet$ Multitask principal-agent problems}

Beyond the single-task moral hazard problem, some multi-task problem (see Holmstr\"{o}m and Milgrom~\cite{HM91} for the discussion of multi-task agency problem) can also be addressed by our approach. Consider a similar model of Bond and Gomes \cite{BG09}, a principal employs an agent to work on $m$ symmetric projects. Each project can result in either success or failure, and the total value of the output is $f_j$, a function of the number of successful projects, $j=0,\dots,m$. The agent's actions determine the probability distribution of success over these tasks: $p=(p_1,\dots,p_m)\in[0,1]^m$. We can view $p$ itself as the agent's action and write $c(p)$ for the agent's cost of action. The limited attention constraint $\sum_{i=1}^mp_i\leq\bar{p}$ is assumed for the agent.

Since the agent's effort allocation is not observable but the number of successful tasks is, the principal rewards the agent based on the number of successes using an incentive scheme $w=(w_j)_{j=0}^m$ where $w_j$ is the agent's reward when there are $j$ successful tasks. Let $A_{m,j}(p)$ be the probability of $j$ successes in the $m$ tasks under action $p$. Then the expected utility functions of the principal and agent are 
\begin{align*}
V(w,p)=\sum_{j=0}^m v(f_j-w_j)A_{m,j}(p),\quad
\mathrm{and}\quad U(w,p)=\sum_{j=0}^m u(w_j)A_{m,j}(p)-c(p),
\end{align*}
respectively, where $v(\cdot)$ and $u(\cdot)$ are the principal's and the agent's preferences over wealth, respectively. In summary, the principal's optimization problem is to select an incentive scheme $w$ and instructions $p$ for the agent so as to maximize his utility subject to the constraints,
\begin{equation*}
\begin{aligned}
\max_{w,p} &\ V(w,p),\\
\mathrm{s.t.}\ &\ p\in\mathrm{arg}\max_{p'\in[0,1]^m} \Big\{U(w,p'): \sum_{i=1}^mp_i'\leq\bar{p}\Big\},\\
&\ U(w,p)\geq\underline{U},\\
&\ \bar{w}\geq w_j\geq0,\ \forall j=0,\dots,m.
\end{aligned}
\end{equation*}
The first constraint is the incentive compatibility constraint, the second is the individual rationality constraint (or participation constraint), and the third is the minimum payment constraint and a budget constraint.

Note that the sequence of probabilities $A_{m,j}(p)$ is a P\'{o}lya frequency sequence in $j$ and it is shown in \cite{BG09} that the sequence $A_{m,j}(p)$ satisfies the monotone likelihood ratio property (MLRP). When the agent is risk neutral, it is a standard result that the principal can maximize his payoff by using an incentive scheme of the form
\begin{align}
w_j=
\left\{
\begin{array}{lll}
0  &\mathrm{if}\ j<k,\\
(1-\mu)\bar{w}  &\mathrm{if}\ j=k,\\
\bar{w}  &\mathrm{if}\ j>k,
\end{array}
\right.
\end{align}
for some integer $k$ and $\mu\in[0,1]$, see, e.g., Proposition 1 in \cite{BG09}. We denote this incentive scheme briefly by $w^{k,\mu}$. For each integer $k\leq m$, after relabeling $x=\mu$ and $y=p$, the principal's contracting problem reduces to \eqref{BLPP} where the upper variable $x$ is one dimensional.

\section{Main generic results and techniques}\label{sec5}

In this section we will first propose {sufficient conditions for the partial calmness conditions in Definitions~\ref{Defn1.2}-\ref{ourcalmness} to hold}. Then  motivated by the  applications discussed in Section \ref{application}, we prove that these conditions are generic for BLPPs under the following problem setting. 

\begin{center}
	\fbox
	{
		\shortstack[l]
		{
			\textbf{Problem Setting:}\\ 
			(1) The real valued mappings $f,g$ are three times continuously differentiable.\\	
			(2) The upper variable $x$ is one dimensional ($n=1$).
		}
	}
\end{center}

Let $J_0(\bar{x},\bar{y}):=\{j: g_j(\bar{x},\bar{y})=0\}$ be the active index set for $\bar{y}\in Y(\bar{x})$. We say that $\bar{y}\in Y(\bar{x})$ is a g.c.  point,  a FJ-point or a KKT point if
\begin{align}
\exists\,(u_0,u)\not =0  \mbox{ s.t. }0&=u_0\nabla_y f(\bar{x},\bar{y})+ \sum_{j\in J_0(\bar{x},\bar{y})}u_j \nabla_y g_j(\bar{x},\bar{y}),  \label{multiplierv}\tag{g.c.}\\
\exists\,(u_0,u)\not =0	 \mbox{ s.t. }	0&=u_0 \nabla_y f(\bar{x},\bar{y})+ \sum_{j\in J_0(\bar{x},\bar{y})} u_j\nabla_y g_j(\bar{x},\bar{y}),  \mbox{ and}\ (u_0, u)\geq 0, \label{FJ} \tag{FJ}\\
\exists\, u \geq 0  \mbox{ s.t. }		0&=\nabla_y f(\bar{x},\bar{y})+ \sum_{j\in J_0(\bar{x},\bar{y})} u_j \nabla_y g_j(\bar{x},\bar{y}), \label{multiplier}\tag{KKT}
\end{align}
respectively. We call $(u_0, u)$ satisfying the g.c. condition a Lagrange vector and the one satisfying KKT condition a (Lagrange) multiplier. Note that  a Lagrange vector is not assumed to be nonnegative.  For each $(x,y)$, denote  the set of FJ-multipliers  at $y$ for problem \ref{Lx} by
\begin{equation}\label{FJmultiplier}
M_{\rm FJ}(x,y):=\left\{(u_0,u)\in\mathbb{R}^{1+p}:
\begin{aligned}
&\ u_0\nabla_{y} f(x,y)+u \nabla_{y}g (x,y)=0, \\
&\ (u_0,u)\geq 0,\  \|(u_0,u)\|_1=1, \ 
u^T g(x,y)=0
\end{aligned}
\right\}.
\end{equation} 
Similarly, we denote the set of KKT-multipliers for \ref{Lx} at $y$ as follows: 
\begin{equation}\label{KKTmultiplier}
M_{\rm KKT}(x,y):=\left\{u:
\nabla_{y} f(x,y)+u \nabla_{y}g (x,y)=0,  \ u\geq 0,\   
u^T g(x,y)=0
\right\}.
\end{equation}

Define 
\begin{align*}
\Sigma_{f}:=&\big\{(x,y)\in\mathbb{R}^{n+m}: y\in Y(x)\big\},\\
\Sigma_{gc} :=&\big\{(x,y)\in\mathbb{R}^{n+m}: y\mathrm{\ is\ a\  g.c.\ point\ for}\ L(x)  \big\},\\
\Sigma_{FJ} :=&\big\{(x,y)\in\mathbb{R}^{n+m}: y\mathrm{\ is\ a\  FJ \ point\ for}\ L(x) \big\},\\
\Sigma_{KKT} :=&\big\{(x,y)\in\mathbb{R}^{n+m}: y\mathrm{\ is\ a\  KKT \ point\ for}\ L(x) \big\},\\
\Sigma_{B} :=&\big\{(x,y)\in\mathbb{R}^{n+m}: y\text{\ is\ a\ B-point\ for}\ L(x) \big\},\\
M:=&\big\{(x,y)\in\mathbb{R}^{n+m}: y\in S(x)\big\}.
\end{align*}
We can show that  
\begin{equation}
M\subseteq \Sigma_B \subseteq 
\Sigma_{FJ} \subseteq \Sigma_{gc} \subseteq \Sigma_{f}, \quad \Sigma_{KKT} 
\subseteq \Sigma_{B},\ \mathrm{and}\ 
\Sigma_B \stackrel{MSCQ}{=} \Sigma_{KKT}.
\label{MFCQinc}
\end{equation}
Indeed, as explained in Section \ref{sec2},   $M\subseteq \Sigma_{B}$ follows by the necessary optimality condition. Now we prove that the  condition $ 0\in \nabla_y f(x,y)+\widehat{N}_{Y(x)}(y)$ implies  FJ condition. Indeed, if MFCQ holds, then it implies  KKT condition (cf.  \cite[Exercise 2.53]{M18} or \cite[Theorem 4.1]{HJO}). Thus there exists a normal multiplier in  FJ condition. Otherwise, by the Motzkin's transposition theorem \cite{M51},  FJ condition with an abnormal multiplier holds. Hence 
\begin{equation*}
0\in \nabla_y f(x,y)+\widehat{N}_{Y(x)}(y) \Longrightarrow  y \mbox{ is\ a\ FJ\  point\ of\ }L(x).
\end{equation*}
To show that $\Sigma_{KKT}\subseteq \Sigma_B$, it suffices to make the following observation. First,  KKT condition for $L(\bar{x})$ at $\bar{y}$ can be rewritten as $-\nabla_{y} f(\bar{x},\bar{y})\in L_{Y(\bar{x})}(\bar{y})^\circ$, where $ L_{Y(\bar{x})}(\bar{y}):=\big\{ d\in\mathbb{R}^m: \langle \nabla_y g_j(\bar{x},\bar{y}),d\rangle\leq0,\ \forall\,j\in J_0(\bar{x},\bar{y})\big\}$ is the linearized cone of $Y(\bar{x})$ at $\bar{y}$. Note that the linearized cone always contains the tangent cone, i.e., $L_{Y(\bar{x})}(\bar{y})\supseteq T_{Y(\bar{x})}(\bar{y})$. Hence $L_{Y(\bar{x})}(\bar{y})^\circ\subseteq  T_{Y(\bar{x})}(\bar{y})^\circ=\widehat{N}_{Y(\bar{x})}(\bar{y})$. This implies that $\Sigma_{KKT}\subseteq \Sigma_B$. The reverse inclusion $\Sigma_{B} \subseteq \Sigma_{KKT}$ usually requires some CQ. Indeed, by \cite[Proposition 4]{GY17}, we have $\Sigma_{B}=\Sigma_{KKT}$ under MSCQ at each point of $(x,y)\in \Sigma_f$.

\subsection{Our main results} 

We propose sufficient conditions to  guarantee the partial calmness conditions in Definitions~\ref{Defn1.2}-\ref{ourcalmness} and show that both of them are generic in Problem Setting. Using the terminology initiated in \cite{LYZ08}, in fact we can define the PEB condition based on other stationary (or feasible) conditions as below. 
\begin{definition}[\bf Partial error bound]\label{vfcq}
	We say that the partial error bound (PEB) over  B,  KKT, FJ, cp, feasible (f)-condition is satisfied at a feasible solution  $(\bar{x},\bar{y})$ of \eqref{BLPP}  if the  partially perturbed  mapping 
	\begin{align*}
	K_{B}(v):=\big\{(x,y): &\ f(x,y)-V(x)\leq v,\, 0\in \nabla_y f(x,y)+\widehat{N}_{Y(x)}(y)\big\},
	\\
	K_{KKT}(v):=\big\{(x,y): &\ f(x,y)-V(x)\leq v,\, \mbox{ y is a KKT point of } L(x) \big\},\\
	K_{FJ}(v):=\big\{(x,y): &\ f(x,y)-V(x)\leq v,\, \mbox{ y is a FJ point of } L(x)\big\},\\
	K_{cp}(v):=\big\{(x,y): &\ f(x,y)-V(x)\leq v,\, \mbox{ y is a g.c.  point of } L(x)\big\},\\
	K_{f}(v):=\big\{(x,y): &\ f(x,y)-V(x)\leq v,\, y\in Y(x)\big\},
	\end{align*}
	is calm at $(0,\bar{x},\bar{y})$, respectively. That is, there exist neighbourhoods $\mathcal{V}$ of $0$ and $\mathcal{U}$ of $(\bar{x},\bar{y})$, and a constant $L\geq 0$ such that 
	\begin{align*}
	\mathrm{dist}_{K(0)}(x,y) \leq L|v|,\ \mathrm{for\ all}\ (x,y)\in\mathcal{U}\cap K(v)\ \mathrm{and\ for\ all\ }v\in\mathcal{V},
	\end{align*} 
	where $K(v):=K_B(v),K_{KKT}(v), K_{FJ}(v), K_{cp}(v), K_f(v)$, respectively.
\end{definition}
Here when we say $y$ is a KKT point, we mean that  KKT condition holds at $y$. Note that  PEB over $f$ is the so-called value function constraint qualification introduced by Henrion and Surowiec \cite{HS11} for BLPPs with convex lower level problems.

\begin{definition}[\bf Uniform weak sharp minimum]
	We say that the uniform weak sharp minimum (UWSM) over  B,  KKT, FJ, cp, feasible (f)-condition is satisfied at a feasible solution  $(\bar{x},\bar{y})$ of \eqref{BLPP}  if there exists a neighbourhood $\mathcal{U}$ of $(\bar{x},\bar{y})$, and a constant $L\geq 0$ such that
	\begin{equation*}
	\mathrm{dist}_{S(x)}(y)\leq L(f(x,y)-V(x)),\ \mathrm{for\ all}\ (x,y)\in\mathcal{U}\cap \Sigma,
	\end{equation*} 
	where $\Sigma=\Sigma_B, \Sigma_{KKT}, \Sigma_{FJ}, \Sigma_{cp}, \Sigma_f$, respectively.
\end{definition}
Note that UWSM over $f$  reduces to the UWSM introduced by Ye and Zhu in \cite{YZ95}. By definition, it is easy to see that the PEB condition is weaker than the UWSM. Furthermore, as pointed out  by Henrion and Surowiec in \cite[Example 3.9]{HS11}, the PEB condition may be strictly weaker than the UWSM.

We now prove that under  PEB condition, the corresponding combined program is partially calm at $(\bar{x},\bar{y})$ as stated in the following theorem. For the case $\Sigma=\Sigma_f$, the result was proved in \cite[Proposition 3.11]{HS11}.
\begin{theorem}\label{partiale}
	Let $(\bar x,\bar y)$ be a local optimal solution of \eqref{BLPP}. If PEB over B, KKT, FJ,  cp, f-condition is satisfied at 
	$(\bar x, \bar y)$, then the corresponding combined program is partially calm at $(\bar x,\bar y)$, respectively. That is,
	there is a constant $\mu\geq 0$ such that $(\bar x,\bar y)$ is a local optimal solution of the partially penalized problem 
	\begin{equation*}
	\begin{aligned}
	\min_{x,y}\ & F(x,y)+\mu\big(f(x,y)-V(x)\big) \\
	\mathrm{s.t.\ \,}
	&(x,y) \in \Sigma, \ G(x,y)\leq0, 
	\end{aligned}
	\end{equation*} 
	where $\Sigma=\Sigma_B,  \Sigma_{KKT}, \Sigma_{FJ}, \Sigma_{cp}, \Sigma_f$, respectively.
\end{theorem}
\begin{proof}
	Since $(\bar x,\bar y)$ is a local solution of \eqref{BLPP} and $F$ is continuously differentiable, without loss of generality, we assume that $F$ is Lipschitz of rank $L_f>0$  on $S:=\Omega\cap \mathcal{U}$, where $\Omega:=\{(x,y): G(x,y)\leq 0\}$ and $\mathcal{U}$ is the {(bounded)} neighbourhood of $(\bar x,\bar y)$ defined in the PEB condition,  and $F(x,y)$ attains a minimum over $C:=M \cap \Omega \cap \mathcal{U}$ at $(\bar x,\bar y)$. By virtue of the exact penalty principle of Clarke (\cite[Proposition 2.4.3]{Clarke}), $F(x,y)+L_f \mathrm{dist}_C(x,y)  $  attains a minimum over $S$. But for all $(x,y)\in \Omega\cap \mathcal{U}$, $\mathrm{dist}_C(x,y)=\mathrm{dist}_{K(0)}(x,y)$. Hence we have
	\begin{equation*}
	\begin{aligned}
	F(\bar x,\bar y) & \leq  F(x,y)+L_f \mathrm{dist}_C(x,y) \quad &\forall (x,y) \in \Omega \cap \mathcal{U}
	\\
	&= F(x,y)+L_f \mathrm{dist}_{K(0)}(x,y)\quad &\forall (x,y) \in \Omega \cap \mathcal{U}\\
	&\leq  F(x,y)+L_f L|v| \quad &\forall (x,y) \in \Omega \cap  \mathcal{U} \cap K(v) \mbox{ and } \forall v\in \mathcal{V}\\
	&=  F(x,y)+L_f L(f(x,y)-V(x))\quad &\forall (x,y) \in \Omega \cap \Sigma \cap  \mathcal{U}  ,
	\end{aligned} 
	\end{equation*}
	where the last equality follows from the fact that $f(x,y)-V(x)\geq 0$ for $y\in Y(x)$. The desired result follows from taking $\mu:=L_fL$.
\end{proof}

By (\ref{MFCQinc}), it is easy to see that 
$$\mbox{PEB over } f \Longrightarrow \mbox{PEB over  g.c. condition}  \Longrightarrow \mbox{PEB over }  FJ  
\Longrightarrow \mbox{PEB over } B. $$ 
If MSCQ holds, by virtue of (\ref{MFCQinc}) again, we have
$$\mbox{PEB over } KKT  \Longleftrightarrow \mbox{PEB over } B.$$
Hence, for the partial calmness condition, the PEB over $B$ is the weakest sufficient condition while PEB over $f$ is the strongest sufficient condition. 

Note that problem \eqref{CP4}  is not a practical problem to solve. However we can show the equivalence  between the partial calmness conditions {for \eqref{CPFJ} and \eqref{CP4}, and} use it to study the generic property of the partial calmness condition for \eqref{CPFJ}.

\begin{theorem}[\bf Equivalence of partial calmness  for \eqref{CPFJ} and \eqref{CP4}]\label{thmequiv}
	Let $\bar y\in S(\bar x)$. Then problem \eqref{CP4} is partially calm at $(\bar x,\bar y)$ if and only if  
	 problem \eqref{CPFJ} is partially calm at $(\bar x,\bar y,\bar u_0, \bar u)$ for each FJ-multiplier $(\bar u_0, \bar u)\in M_{\rm FJ}(\bar{x},\bar{y})$. 
\end{theorem}
\begin{proof}
	Suppose that  \eqref{CP4} is partially calm at $(\bar x,\bar y)$. Then there exist a calmness modulus $\bar \mu\geq0$ and a neighborhood $U(\bar{x},\bar{y})$ of $(\bar{x},\bar{y})$ such that 
	\begin{equation}\label{key1}
	F(\bar{x},\bar{y})\leq F(x,y)+\bar \mu\big(f(x,y)-V(x)\big), \ \forall\,(x,y)\in \Sigma_{FJ}\cap U(\bar{x},\bar{y}).
	\end{equation}
	Let  $(x',y', u_0',  u')$ be a feasible solution for problem \eqref{CP2} with $(x',y')$ lying in the neighborhood  $U(\bar{x},\bar{y})$.
	Then $(x',y')\in \Sigma_{FJ}\cap U(\bar{x},\bar{y})$. Hence (\ref{key1}) implies that problem \eqref{CPFJ} is partially calm at $(\bar x,\bar y,\bar u_0,\bar u)$ for any FJ-multiplier $(\bar u_0, \bar u)\in M_{\rm FJ}(\bar x,\bar y)$ with the same calmness modulus $\mu=\bar \mu$.
	
	Conversely,  suppose that problem \eqref{CPFJ} is partially calm at $(\bar x,\bar y,\bar u_0, \bar u)$ for each FJ-multiplier $(\bar u_0, \bar u)\in M_{\rm FJ}(\bar x,\bar y)$. 
	We first show that there exists $\bar{\mu}\geq0$ such that  \eqref{CPFJ} is partially calm at $(\bar x,\bar y,\bar u_0, \bar u)$ for all FJ-multipliers $(\bar u_0, \bar u)\in M_{\rm FJ}(\bar x,\bar y)$ with the same $\mu=\bar{\mu}$.
	 To the contrary,  suppose  that there exists  $(\bar{u}_0^k,\bar{u}^k)\in M_{\rm FJ}(\bar x,\bar y)$ such that problem \eqref{CPFJ} is partially calm at $(\bar x,\bar y,\bar u_0^k, \bar u^k)$ and as $k\rightarrow\infty$, 
	\begin{equation*}
	\mu^k:=\inf\Big\{\mu\geq0:\text{problem \eqref{CPFJ} is partially calm at $(\bar x,\bar y,\bar u_0^k, \bar u^k)$ with $\mu$}\Big\}\rightarrow \infty
	\end{equation*}
 Since $M_{\rm FJ}(\bar x,\bar y)$ is compact, up to a subsequence, the sequence $(\bar{u}_0^k,\bar{u}^k)$ has a limit $(\bar{u}_0^*,\bar{u}^*)$ in $M_{\rm FJ}(\bar x,\bar y)$. By  assumption, 
	\eqref{CPFJ} is partially calm at $(\bar x,\bar y,\bar{u}_0^*,\bar{u}^*)$. Hence, 
	there is $\mu^*\geq0$ and a neighborhood $U(\bar{x},\bar{y},\bar{u}_0^*,\bar{u}^*)$ of $(\bar{x},\bar{y},\bar{u}_0^*,\bar{u}^*)$ such that     
	\begin{equation}
	F(\bar{x},\bar{y})\leq F(x,y)+\mu^*\big(f(x,y)-V(x)\big), \ \forall\,(x,y,u_0,u)\in \mathcal{F}_{\mathrm{CPFJ}}\cap U(\bar{x},\bar{y},\bar{u}_0^*,\bar{u}^*),\label{CPFJcalmness}
	\end{equation}
	where $\mathcal{F}_{\mathrm{CPFJ}}$ is the feasible region of \eqref{CP2}.
	Since $(\bar{u}_0^k,\bar{u}^k)\rightarrow(\bar{u}_0^*,\bar{u}^*)$ in $M_{\rm FJ}(\bar x,\bar y)$,
	 for sufficiently large $k$, 
	we can take a neighborhood $U(\bar{x},\bar{y},\bar{u}_0^k,\bar{u}^k)$ of $(\bar{x},\bar{y},\bar{u}_0^k,\bar{u}^k)$ such that $U(\bar{x},\bar{y},\bar{u}_0^k,\bar{u}^k)\subseteq U(\bar{x},\bar{y},\bar{u}_0^*,\bar{u}^*)$. Hence (\ref{CPFJcalmness}) holds with $U(\bar{x},\bar{y},\bar{u}_0^*,\bar{u}^*)$ replaced with $U(\bar{x},\bar{y},\bar{u}_0^k,\bar{u}^k)$.
	This implies that problem \eqref{CPFJ} is partially calm at $(\bar x,\bar y,\bar{u}_0^k,\bar{u}^k)$ with $\mu=\mu^*$ for all $k$ sufficiently large. Thus $\mu^k\leq \mu^*$ for sufficiently large $k$, a contradiction with $\mu^k\rightarrow\infty$ as $k\rightarrow\infty$.
	
	Now we show that  \eqref{CP4} is partially calm at $(\bar x,\bar y)$ with $\mu=\bar{\mu}$. Suppose not, then there exists $(x^k,y^k)\in\Sigma_{FJ}$ satisfying $(x^k,y^k)\rightarrow(\bar{x},\bar{y})$ and 
	\begin{equation}\label{nonmin}
	F(\bar{x},\bar{y})> F(x^k,y^k)+\bar{\mu}\big(f(x^k,y^k)-V(x^k)\big),\quad\forall\, k.
	\end{equation}
	Take $(u_0^k,u^k)\in M_{\rm FJ}(x^k,y^k)$. Since
	$
	M_{\rm FJ}(x^k,y^k)\subseteq\big\{(u_0,u): \|(u_0,u)\|_1=1\big\},
$
	the sequence $(u_0^k,u^k)\rightarrow(\bar{u}_0,\bar{u})$, up to a subsequence. 
	It follows that $(\bar{u}_0,\bar{u})\in M_{\rm FJ}(\bar x,\bar y)$  due to the continuous differentiability of $ f(x,y) $ and $ g(x,y)$.
	 Hence,
	   $(x^k,y^k,u_0^k,u^k)\rightarrow(\bar{x},\bar{y},\bar{u}_0,\bar{u})$ and \eqref{nonmin} holds. Therefore,  $(\bar{x},\bar{y},\bar{u}_0,\bar{u})$ is not a local minimum of \eqref{CP2} with $\mu=\bar{\mu}$, a contradiction.	
\end{proof}

The results in Theorem~\ref{thmequiv} can be improved if an extra assumption is imposed as below.
\begin{corollary}
	Let $\bar y\in S(\bar x)$. Assume that the constant rank constraint qualification (CRCQ)  for the lower level problem holds at $(\bar x,\bar y)$, i.e., there exists an open neighborhood $U(\bar{x},\bar{y})$ of $(\bar{x},\bar{y})$ such that for each index set $I\subseteq J_0(\bar{x},\bar{y})$ the family of gradients $\{\nabla_{y} g_j(x,y): j\in I\}$ has the same rank on $U(\bar{x},\bar{y})$.
	 Then   the combined program \eqref{CP4}  is partially calm at $(\bar x,\bar y)$ if problem \eqref{CPFJ} is partially calm at $(\bar{x},\bar{y},\bar{u}_0,\bar{u})$ for all $(\bar{u}_0,\bar{u})$ which are  vertices of $ M_{\rm FJ}(\bar x,\bar y)$. 
\end{corollary}
\begin{proof}
	The proof is similar to the argument of the proof of   \cite[Corollary 3.3]{DD12}. Indeed,  FJ-multipliers $(u_0^k,u^k)$ in the proof of Theorem~\ref{thmequiv} can without assumptions be taken as vertices of $M_{\rm FJ}(x^k,y^k)$. Then CRCQ implies that the sequence $(u_0^k,u^k)$ converges to a vertex of $M_{\rm FJ}(\bar x,\bar y)$. Hence, the desired result follows.
\end{proof}

Through the equivalence between the partial calmness for \eqref{CPFJ} and \eqref{CP4}, we can study the generic property of the partial calmness condition for \eqref{CPFJ}.

We can now formulate our main results. 
\begin{theorem}[\bf Generic property]\label{thm1}
	Under Problem Setting, we have  
	\begin{enumerate}[(1)]
		\item[{\rm (1)}] Let $\mathcal{P}$ denote the set of functions $(f,g)\in [C^3(\mathbb{R}^{n+m})]^{1+p}\cap\mathcal{O}$ such that the partial error bound over  {FJ-condition} defined in Definition \ref{vfcq} holds at all points $(\bar{x},\bar{y})\in M$. Then $\mathcal{P}$ is $C^3_s$-open and $C^3_s$-dense in $[C^3(\mathbb{R}^{n+m})]^{1+p}\cap\mathcal{O}$. 
		Consequently, the partial error bound over {FJ}-condition is generic.
		
		\item[{\rm (2)}] Let $\mathcal{Q}$ denote the set of functions $(f,g)\in [C^3(\mathbb{R}^{n+m})]^{1+p}\cap\mathcal{O}$
		such that the partial calmness for \eqref{CP4} holds at all local optimal solutions of \eqref{BLPP}. Then $\mathcal{Q}$ is $C^3_s$-open and $C^3_s$-dense in $[C^3(\mathbb{R}^{n+m})]^{1+p}\cap\mathcal{O}$. Consequently, the partial calmness condition for \eqref{CPB} and  \eqref{CP4} are generic.
	
		\item[{\rm (3)}] Let $\mathcal{R}$ denote the set of functions $(f,g)\in [C^3(\mathbb{R}^{n+m})]^{1+p}\cap\mathcal{O}$
		such that the partial calmness for \eqref{CPFJ} holds at all local optimal solutions  of \eqref{CPFJ}. Then $\mathcal{R}$ is $C^3_s$-open and $C^3_s$-dense in $[C^3(\mathbb{R}^{n+m})]^{1+p}\cap\mathcal{O}$. Consequently, the partial calmness condition for  \eqref{CPFJ} is generic.
	\end{enumerate}
	Here $\mathcal{O}$ stands for the set of $(f,g)$ such that 
	\begin{equation}\label{JS24}
	B_{f,g}(\bar{x},c)\mathrm{\ is\ compact\ for\ all\ }(\bar{x},c)\in\mathbb{R}^n\times\mathbb{R},
	\end{equation}
	where 
	$
	B_{f,g}(\bar{x},c):=\big\{(x,y): \|x-\bar{x}\|\leq 1,\, f(x,y)\leq c, \, g(x,y)\leq0\big\}.
	$
\end{theorem}

The proof  for part (1) will be given in Section \ref{sec6}. Assuming Part (1) is true, part (2) follows from the fact that the partial calmness condition for \eqref{CPB} is weaker than the one for \eqref{CP4}. {Part (3)} follows from the equivalence between  the partial calmness  for \eqref{CPFJ} and \eqref{CP4} in Theorem 
\ref{thmequiv}.

\begin{remark} It is worth pointing out that the generic analysis can not be done with the set $\Sigma_f$. Hence the value function CQ of  Henrion and Surowiec \cite{HS11} is not generic and neither is the partial calmness for \eqref{VP} (i.e., the partial calmness over $f$). 
	
The condition \eqref{JS24} (also in \cite{JS12}) implies the inf-compactness condition: for all $\bar{x}\in\mathbb{R}^n$ there exist $\alpha>0$, $\delta>0$ and a bounded set $C$ such that $\alpha>V(\bar{x})$ and 
$
	\big\{y: \|x-\bar{x}\|\leq\delta,\, f(x,y)\leq\alpha, \, g(x,y)\leq0\big\}\subseteq C.
$
Since $\bar{y}\in S(\bar{x})$ is a global minimizer for $L(\bar{x})$, the inf-compactness condition is sufficient to avoid asymptotical effects for the feasible set of BLPPs. Note also  that $[C^3(\mathbb{R}^{n+m})]^{1+p}\cap\mathcal{O}$ is $C^3_s$-open.
\end{remark}

\subsection{Techniques and Roadmap of Analysis}\label{Sec4.2}

It is clear that the main goal of Theorem~\ref{thm1} is  to verify that  the partial error bound over FJ-condition defined in Definition \ref{vfcq} holds generically. The main difficulty is to deal with the global constraint $y\in S(x)$ or the value function constraint $f(x,y)-V(x)\leq0$. Hence the generic structure of $M=\{(x,y): y\in S(x)\}$ is important. 

To describe the generic structure of $M$, it is natural to adopt a parametric programming point of view, using the upper level variable as a parameter. This issue has been studied by Jongen and Shikhman \cite{JS12} when the upper level variable is one dimensional. The main result in \cite{JS12} states that--generically--the set $M$ is the union of $C^2$ curves with boundary points and kinks which can be parametrized by means of the upper variable. The appearance of the boundary points and kinks is due to certain degeneracies of optimal solutions for the lower level problem (e.g., the failure of strict complementarity (SC) or LICQ or the second order condition (SOC), as well as the change from local to global solutions).

The main idea in \cite{JS12} is to exploit a generic classification of g.c. points for one dimensional parametric optimization problems in \cite{JJT86}. When $n=1$, it was shown in \cite{JJT86} that for an open and dense subset of problem data (i.e., generically) the set $\Sigma_{gc}$ is pieced together from one-dimensional manifolds. Moreover, for those problem data, the g.c. points for $L(x)$ can be divided into five (characteristic) types.  A point $(\bar{x},\bar{y})\in\Sigma_{gc}$ is of type 1 if it is a nondegenerate critical point for $L(\bar{x})$ which  means that the following conditions hold:
\begin{enumerate}[ND1:]
	\item LICQ is satisfied at $(\bar{x},\bar{y})$ and hence there exists a unique multiplier $\bar u$.
	
	\item SC holds at $(\bar{x},\bar{y})$, i.e., $\bar u_j\not =0,\  j\in J_0(\bar x,\bar y)$.
	
	\item SOC is satisfied, i.e., the matrix 
	\begin{equation} 
	(\nabla_{yy}^2( f+  g^T\bar u)(\bar x,\bar y))|_{T_{\bar{y}}Y(\bar{x})}:= V^T\nabla_{yy}^2 (f+g^T \bar u)(\bar x,\bar y)V\label{matrix}
	\end{equation}  
	is nonsingular where   $V$ is some matrix whose columns form a basis for the tangent space $T_{\bar{y}}Y(\bar{x}):=\{ \xi\in \mathbb{R}^m :  \nabla_y g_j(\bar x,\bar y) ^T \xi =0,\  j\in J_0(\bar x,\bar y)\}.$
\end{enumerate}
Other types originate from degeneracies of the generalized critical points. A degeneracy describes the failure of at least one of the nondegeneracy properties ND1-ND3. For cleanness of the paper, we do not repeat the precise characterization of other types here and introduce them in Table \ref{5types} for brevity. For more details we refer the reader to \cite{JJT86,JS12}.
\begin{table}[h!]
	\begin{center}		
		\scalebox{0.8}{
			\begin{tabular}{cc}
				& \textbf{Characteristic}\\
				\toprule
				\textbf{Type 1} &  LICQ, SC, and SOC are satisfied\\
				\midrule
				\textbf{Type 2} &  LICQ, SOC are satisfied, SC is violated and exactly one of  
				${\bar u}_j, j\in J_0$ vanishes \\
				\midrule
				\textbf{Type 3} &  LICQ, SC are satisfied, and exactly one eigenvalue of the matrix (\ref{matrix}) vanishes\\
				\midrule
				\textbf{Type 4} &  LICQ violated, the rank of $\{ \nabla_y g_j(\bar x,\bar y),\, j\in J_0\} =|J_0|-1<m$, and KKT condition fails\\
				\midrule
				\textbf{Type 5-1} & \makecell[c]{LICQ violated, the rank of $\{ \nabla g_j(\bar x,\bar y),\, j\in J_0\} =|J_0|=m+1$,\\ MFCQ fails but KKT condition holds}\\
				\midrule
				\textbf{Type 5-2} & \makecell[c]{LICQ violated, the rank of $\{ \nabla g_j(\bar x,\bar y),\,j\in J_0\} =|J_0|=m+1$, and MFCQ holds }\\
				\bottomrule
		\end{tabular}}
		\caption{The Five Types in Parametric Optimization $(n=1)$. Here $J_0:=J_0(\bar x,\bar y)$.}
		\label{5types}
	\end{center}
\end{table}

We next briefly review the main result of Jongen and Shikhman \cite{JS12} for BLPPs. To this end, we first recall the definition of simplicity of bilevel problems in \cite{JS12}. 

\begin{definition}[\bf \text{\cite[Definition 4.1]{JS12}}]\label{types}
	A bilevel programming problem (with $n=1$) is called simple at $(\bar{x},\bar{y})\in M$ if one of the following cases occurs: 
	
		Case I: $S(\bar{x})=\{\bar{y}\}$ and $(\bar{x},\bar{y})$ is of Type 1, 2, 4, 5-1 or 5-2. 
		
		Case II: $S(\bar{x})=\{\bar{y}_1,\bar{y}_2\}$ and $(\bar{x},\bar{y}_1)$, $(\bar{x},\bar{y}_2)$ are both of Type 1. In addition,  it holds that
		\begin{equation}\label{alpha}
		\alpha:=\frac{d}{dx}[f(x,y_2(x))-f(x,y_1(x))]\Big|_{x=\bar{x}}
		\neq0,
		\end{equation}
		where $y_1(x)$, $y_2(x)$ are unique local minimizers for $L(x)$ in a neighborhood of $\bar{x}$ with $y_1(\bar{x})=\bar{y}_1$, $y_2(\bar{x})=\bar{y}_2$ according to Type 1. 
\end{definition}

\begin{theorem}[\bf \text{\cite[Theorem 4.1]{JS12}}]\label{thm41}
	
	Let $\mathcal{F}$ be the set of defining functions $(F,f,g)\in C^3(\mathbb{R}^{1+m})\times\big([C^3(\mathbb{R}^{1+m})]^{1+p}\cap\mathcal{O}\big)$ such that the corresponding BLPP is simple at all its feasible points $(\bar{x},\bar{y})\in M$. Then $\mathcal{F}$ is $C^3_s$-open and $C^3_s$-dense in $C^3(\mathbb{R}^{1+m})\times\big([C^3(\mathbb{R}^{1+m})]^{1+p}\cap\mathcal{O}\big)$. 
\end{theorem}

\begin{remark}\label{rk1}
		First, since the main goal of Jongen and Shikhman \cite{JS12} is to describe the generic structure of the bilevel feasible set $M$, they focus on such parts of $\Sigma_{gc}$ which correspond to (local) minimizers, i.e.,
		\begin{equation*}
		\Sigma_{min}:=\big\{(x,y)\in\mathbb{R}^{n+m}: y\mathrm{\ is\ a\  local\ minimizer\ for}\ L(x)  \big\}
		\end{equation*}
		in a neighborhood of $(\bar{x},\bar{y})\in M$. Hence,  the cases in \cite[Definition 4.1]{JS12} are related to $\Sigma_{min}$; cf. \cite[Section 3]{JS12}. 
		
		Our purpose is to study the partial calmness conditions for the combined programs \eqref{CPB} and \eqref{CP4} which involve B- and FJ conditions. Note that  $\Sigma_{min}\subseteq\Sigma_{B}\subseteq\Sigma_{FJ}$. Hence it is crucial whether \cite[Theorem 4.1]{JS12} is still true when the set $\Sigma_{min}$ in Definition \ref{types} is replaced by $\Sigma_{FJ}$. The answer is affirmative. Actually, a stronger result is in fact true. The definition of simplicity, those analysis in \cite[Section 3]{JS12} and then \cite[Theorem 4.1]{JS12} can be generalized to a so-called `Kuhn-Tucker subset' $\Sigma_{KT}$ of $\Sigma_{gc}$. This notation is studied in more details, in particular in connection with the (failure of) MFCQ, cf. \cite[Section 4]{JJT86}. More precisely, the subset $\Sigma_{KT}$ of $\Sigma_{gc}$ is defined to be the closure of nondegenerate critical points (i.e., points of Type 1) with nonnegative Lagrange multipliers. 
		
		Now we show that $(x,y)\in\Sigma_{FJ}$ and $(x,y)$ is of Types 1-5 imply that $(x,y)\in\Sigma_{KT}$. Indeed, point $(x,y)\in\Sigma_{FJ}$ implies that the multipliers are nonnegative. Thus, if further $(x,y)$ is of Type 1, then it is clear that  $(x,y)\in\Sigma_{KT}$. If $(x,y)$ is of Type 2, then $y$ is a nondegenerate critical point for the problem $\widetilde{L}(\bar{x})$, which differs from $L(\bar{x})$ by omitting the inequality constraint with vanishing Lagrange multiplier. Hence, we get a unique $C^2$-map $\widetilde{y}(x)$. But only one of its branches from $\bar{x}$ satisfies the omitting inequality constraint with strictly sign, and then belongs to the set of Type 1 points with positive Lagrange multipliers. Hence $(\bar{x},\bar{y})$ belongs to $\Sigma_{KT}$. Similar arguments apply to the left cases.  Since $\Sigma_{min} \subseteq \Sigma_{FJ}$, Theorem \ref{thm41} is still true when the set $\Sigma_{min}$ in the definition \ref{types} is replaced by $\Sigma_{FJ}$.
\end{remark}

It is worth pointing out that the proof of \cite[Theorem 4.1]{JS12} does not use any property of the upper level objective function. In other words, the generic structures of $\Sigma_{B},\Sigma_{KKT},\Sigma_{FJ},\Sigma_{KT},\Sigma_{gc},\Sigma_{min}$ are properties of the  lower level \eqref{CP4} only, independent of the upper level objective function. Hence, \cite[Theorem 4.1]{JS12} can be stated as follows. 

\begin{theorem}[\bf Modified Theorem 4.1 in \cite{JS12}]\label{newthm41}	
	Let $\mathcal{F}$ denote the set of defining functions $(f,g)\in \big([C^3(\mathbb{R}^{1+m})]^{1+p}\cap\mathcal{O}\big)$ such that the corresponding bilevel programming problem is simple at all its feasible points $(\bar{x},\bar{y})\in M$. Then $\mathcal{F}$ is $C^3_s$-open and $C^3_s$-dense in $\big([C^3(\mathbb{R}^{1+m})]^{1+p}\cap\mathcal{O}\big)$. 
\end{theorem}

\begin{figure}[h!]
	\centering
	
	\includegraphics[scale=0.2]{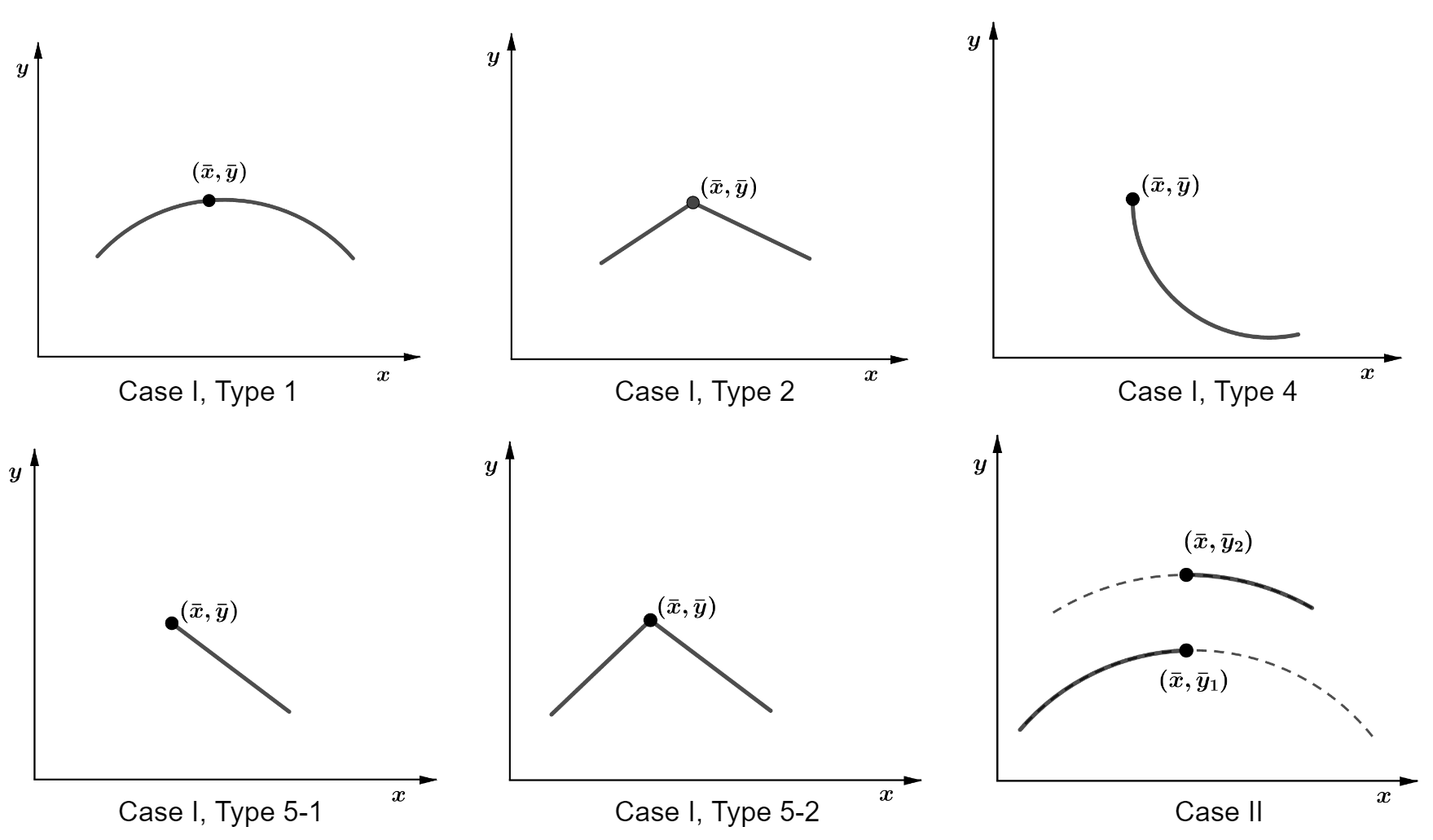}
	\caption{Generic structure of FJ points locally around global minimizers.}
	\label{figtypes}
\end{figure}

For the convenience of the reader, in each case one of the possibilities in Definition \ref{types} is depicted in Figure \ref{figtypes}, which can be seen as an extension of Fig. 1 in \cite{JS12}. 

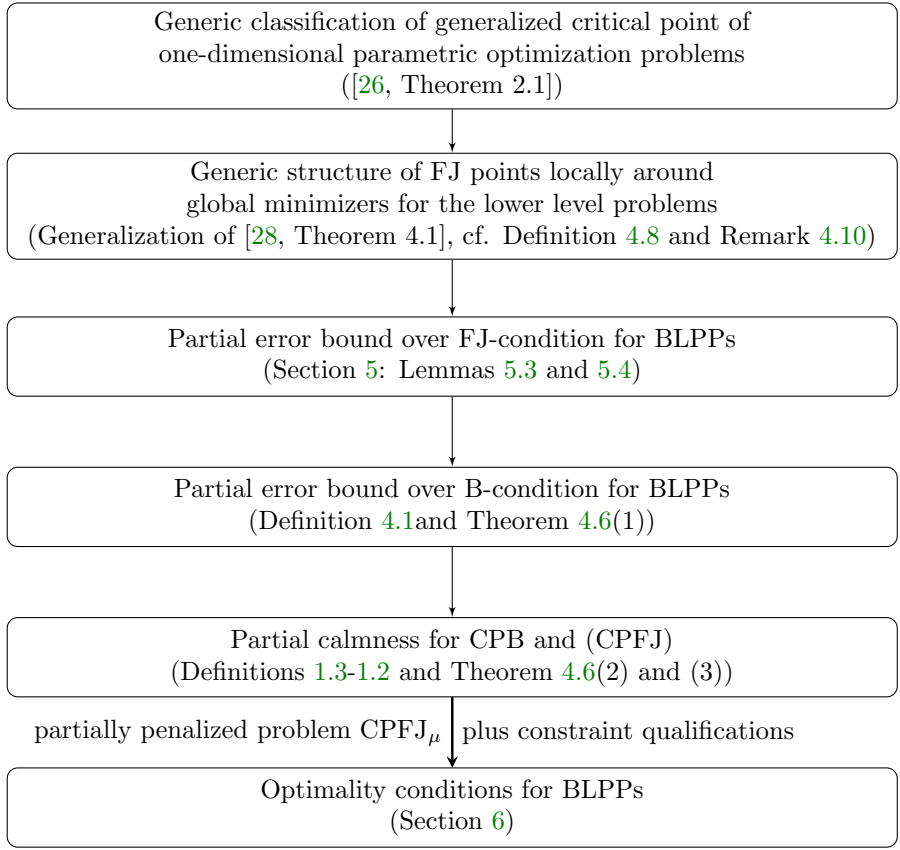
\begin{figure}[h!]
	\centering
	\begin{tikzpicture}[
	minimum size=1cm,
	node distance=2cm and 2.1cm,
	>=stealth,
	bend angle=45,
	auto
	]
	\node [block] (init) {Generic classification of generalized critical point of one-dimensional parametric optimization problems\\
		(\cite[Theorem 2.1]{JJT86})};  
	\node [block, below of= init](ReadA){Generic structure of FJ points locally around global minimizers for the lower level problems\\
		(Generalization of \cite[Theorem 4.1]{JS12},
		cf. Definition \ref{types} and Remark \ref{rk1})};  
	\node [block, below of= ReadA](ReadB){Partial error bound over FJ-condition for BLPPs\\
		(Section \ref{sec6}: Lemmas \ref{lem6.2} and \ref{lem6.3})};  
	\node [block, below of = ReadB](Sum){Partial error bound over B-condition for BLPPs\\
		(Definition \ref{vfcq}and Theorem \ref{thm1}(1))};  
	\node [block, below of = Sum](P){Partial calmness for CPB and (CPFJ)\\
		(Definitions \ref{ourcalmness}-\ref{Defn1.2}  and Theorem \ref{thm1}(2) and (3))};  
	\node [block, below of = P](Out){Optimality conditions for BLPPs\\
		(Section \ref{sec7})};  
	\path [line] (init) -- (ReadA);   
	\path [line] (ReadA) -- (ReadB);   
	\path [line] (ReadB) -- (Sum);   
	\path [line] (Sum) -- (P);   
	\draw [arrow] (P) -- node[anchor=west] {plus  constraint qualifications} (Out);   
	\draw [arrow] (P) -- node[anchor=east] {partially penalized problem CPFJ$_\mu$} (Out);
	\end{tikzpicture}  
	\caption{Roadmap to study the genericity result.}
	\label{roadmap}
\end{figure}

To present our analysis more clearly, let us show the roadmap in Figure \ref{roadmap}. 

\section{Proof of the main result}\label{sec6}

The main aim of this section is to prove Theorem \ref{thm1}(1). By Theorems \ref{partiale} and \ref{newthm41}, and the relationship between the partial error bounds over FJ and B conditions, it suffices to verify partial error bound over FJ-condition holds in Case I and Case II.

To verify partial error bound over FJ-condition, the following lemma is useful. 

	\begin{lemma}\label{lem1}
	Assume that both $f(x,y)$ and $g(x,y)$ are continuous. Let $\bar{y}\in S(\bar{x})$ and $x_i\rightarrow\bar{x}$ as $i\rightarrow\infty$. Suppose that there exist $\bar{y}_i\in Y(x_i)$ such that $\bar{y}_i\rightarrow \bar{y}$. Then $y_0\in S(\bar{x})$ if $y_i\in S(x_i)$ and $y_i\rightarrow y_0$.
\end{lemma}

\begin{proof}
	Let $y_i\in S(x_i)$ and $\lim_{i\rightarrow \infty }y_i=y_0$. Since $y_i\in S(x_i)\subseteq Y(x_i)$, we have $g(x_i,y_i)\leq0$. Thus, by the continuity of $g(x,y)$,   
	$
	g(\bar{x},y_0)=\lim_{i\rightarrow \infty } g(x_i,y_i)\leq0.
	$
	Hence $y_0\in Y(\bar{x})$. Next we show that $f(\bar{x},y_0)\leq f(\bar{x},\bar{y})$. Indeed, since $y_i\in S(x_i)$ and $\bar{y}_i\in Y(x_i)$, we have $f(x_i,y_i)\leq f(x_i,\bar{y}_i)$. Hence, by the continuity of $f(x,y)$, 
	$
	f(\bar{x},y_0)=\lim_{i\rightarrow \infty }f(x_i,y_i)
	\leq\lim_{i\rightarrow \infty }f(x_i,\bar{y}_i)=f(\bar{x},\bar{y}).
	$
	Therefore, it is clear that $y_0\in S(\bar{x})$ since $\bar{y}\in S(\bar{x})$.  
\end{proof}

\begin{remark}
	The above lemma is important enough to be stated separately since we will apply it to verify UWSM over FJ-condition in Case I and PEB over FJ-condition in Case II. Its principal significance is that it allows us to show that some branch of $\Sigma_{FJ}$ become global minimizers under some condition (e.g., $|S(\bar{x})|\leq2$).   
\end{remark}

\subsection{Verification of UWSM over FJ-condition in Case I}
In fact in Case I, we can show the stronger condition UWSM {over FJ-condition} holds and the corresponding constant $L=0$.
\begin{lemma}\label{lem6.2}
	Let \eqref{BLPP} be simple at $(\bar{x},\bar{y})\in M$ in Case I. Then there exists a neighbourhood $\mathcal{U}$ of $(\bar{x},\bar{y})$ such that 
	{ $(x,y)\in \mathcal{U}$ and $y$ is a FJ-point for $L(x)$ imply that $y$ minimizes $L(x)$, i.e., 
		\begin{equation}\label{fjm}
		\Sigma_{FJ}\cap\mathcal{U}=M\cap\mathcal{U}.
		\end{equation}}
	Hence for all $(x,y)\in \Sigma_{FJ}\cap \mathcal{U}$, 
	\begin{equation}
	\mathrm{dist}_{S(x)}(y)=0. 
	\end{equation}
	This implies that the UWSM over FJ condition holds  at $(\bar{x},\bar{y})$ and the corresponding constant $L=0$.
\end{lemma}

\begin{proof}
	The proof falls naturally into four parts.
	
	\textbf{Part 1: Case I, Type 1.}
	In this case, a stronger result is in fact true. We claim that there exists a neighbourhood $\mathcal{U}$ of $(\bar{x},\bar{y})$ such that
	\begin{equation}\label{gcm}
	\Sigma_{gc}\cap\mathcal{U}=M\cap\mathcal{U}.
	\end{equation}
	The proof will be divided into two steps. 
	
	\noindent
	{\bf Step 1}: There exists a neighbourhood $\mathcal{U}:=\mathcal{U}(\bar{x}) \times \mathcal{U}(\bar{y})$ of $(\bar{x},\bar{y})$ such that for some continuous function $y(x)$, 
	$
	\Sigma_{gc}\cap  \mathcal{U}= \big\{(x,y(x)):x\in \mathcal{U}(\bar{x})\big\}.
	$

	This result is standard, cf. \cite[Section 3.1]{JJT86}. Indeed, it is well-known that conditions ND1-ND3 allow us to apply the implicit function theorem and obtain a unique $C^2$-mapping $y(x)$ in an open neighborhood of $\bar{x}$ such that we can parametrize the set $\Sigma_{gc}$ by means of a unique $C^2$-mapping $x\mapsto (x,y(x))$ in an open neighbourhood of $(\bar{x},\bar{y})$. 
	As a consequence, the set $Y(x)\neq\emptyset$ for all $x\in\mathcal{U}(\bar{x})$. It then follows from the definition of $\mathcal{O}$, the solution set $S(x)\neq\emptyset$ for all $x\in\mathcal{U}(\bar{x})$. 
	
	\noindent 	{\bf Step 2}: There exists a small neighbourhood $\mathcal{U}(\bar{x})$ of $\bar{x}$ such that $y(x)$ minimizes $L(x)$ for all $x\in\mathcal{U}(\bar{x})$. 
	
	Suppose not, then there exists $x_i\rightarrow\bar{x}$ such that $y(x_i)\notin S(x_i)$. Taking $\bar{y}_i=y(x_i)$, by the continuity of $y(x)$, we have $\bar{y}_i\rightarrow\bar{y}$. Taking $y_i\in S(x_i)$ and $y_i\rightarrow y_0$, by Lemma \ref{lem1}, we have $y_0\in S(\bar{x})$. Recall that $S(\bar{x})=\{\bar{y}\}$ in Case I. Hence $y_0=\bar{y}$ and then $y_i\rightarrow\bar{y}$. By the fact that optimal solution $y_i$ must be a g.c. point for $L(x_i)$. Hence, when $i$ is sufficient large, by Step 1, it is immediate that $y_i=y(x_i)$ and then $y(x_i)\in S(x_i)$ since $y_i\in S(x_i)$. This contradicts our assumption that $y(x_i)\notin S(x_i)$.   
	
	\textbf{Part 2: Case I, Type 2.}
	In this case, LICQ  is satisfied and the only degeneracy is the vanishing of exactly one Lagrange multiplier corresponding to an active inequality constraint, indexed by $j_0$. Then $\bar{y}$ is a nondegenerate critical point for the problems $\widetilde{L}(\bar{x})$ and $\widehat{L}(\bar{x})$, where $\widetilde{L}(\bar{x})$ differs from $L(\bar{x})$ by omitting the inequality constraint $g_{j_0}$ and whereas in  $\widehat{L}(\bar{x})$ the constraint $g_{j_0}$ is regarded as an equality constraint. Hence, the set $\Sigma_{gc}$ around $(\bar{x},\bar{y})$ consists of the feasible part of two curves, cf. \cite[Section 3.2]{JJT86}.
	
	Let the KT subset $\Sigma_{KT}$ of $\Sigma_{gc}$ be the closure of the set of Type 1 points with nonnegative Lagrange multipliers. Because $\bar{y}\in S(\bar{x})$ implies that $\bar{y}$ is a local minimizer for $L(\bar{x})$, similar to those analysis in \cite[Section 3.2]{JS12}, by Morse relations, we can parameterize $\Sigma_{KT}\cap \mathcal{U}$ by the $x$-variable. Hence Step 1 in Part 1 is still true for some continuous function $y(x)$ when $\Sigma_{gc}$ is replaced by $\Sigma_{KT}$. Note that the continuity of $y(x)$ is enough in the proof of Step 2 in Part 1. Therefore, by the same argument in Part 1, we can prove 
	\begin{equation}\label{ktm}
	\Sigma_{KT}\cap\mathcal{U}=M\cap\mathcal{U}.
	\end{equation}

	Now we show that there exists a neighbourhood $\mathcal{U}$ of $(\bar{x},\bar{y})$ such that  $(x,y)\in \mathcal{U}$ and $y$ is a FJ-point for $L(x)$ imply that $y$ minimizes $L(x)$. By the argument in the proof of Step 2 in Part 1 and \cite[Theorem 2.1]{JJT86}, since an optimal solutin must be a FJ point, it suffices to show that $(x,y)\in\Sigma_{FJ}$ and $(x,y)$ is of Types 1-5 imply that $(x,y)\in\Sigma_{KT}$. In fact, this result has been shown in Remark \ref{rk1}.
	
	\textbf{Part 3: Case I, Type 4.}
	In this case, LICQ is violated and the number of active constraints is less than $m+1$. Moreover, $(\bar{x},\bar{y})$ is no longer a KKT-point. But in an open neighborhood of $(\bar{x},\bar{y})$, apart from $(\bar{x},\bar{y})$, the points of $\Sigma_{gc}$ are nondegenerate critical points. Let the parameter $x$ viewed as a function on $\Sigma_{gc}$, the set $\Sigma_{gc}$ can be locally approximated by means of a parabola, cf. \cite[Section 3.4]{JJT86}. 
	
	Because $\bar{y}$ is a local minimizer for $L(\bar{x})$ when $\bar{y}\in S(\bar{x})$, similar to those analysis in \cite[Section 3.4]{JS12}, it can be shown that only one branch of $\Sigma_{gc}$ starting from $\bar{y}$ left in $\Sigma_{KT}$. Hence we can prove \eqref{ktm} by the same argument in Part 2 and then obtain the desired result.
	
	\textbf{Part 4: Case I, Type 5-1 or 5-2.}
	In this case, LICQ is violated, but in contract to Type 4, the number of active constraints equals $m+1$. Then, by the characteristic conditions in Type 5 (cf. \cite[Section 3.5]{JJT86}), around $(\bar{x},\bar{y})$, the set $\Sigma_{gc}$ consists of exactly $m+1$ (half) curves, each of them either emanating from $(\bar{x},\bar{y})$ or ending at $(\bar{x},\bar{y})$. Moreover, there exists an open neighbourhood $\mathcal{U}$ of $(\bar{x},\bar{y})$ such that either only one branch of $\Sigma_{KT}\cap \mathcal{U}$ left (referred to this case as Type 5-1) or one emanating from $(\bar{x},\bar{y})$ with another one ending at $(\bar{x},\bar{y})$ in $\Sigma_{KT}$ (referred to this case as Type 5-2), cf. \cite[Section 3.5]{JJT86}. Furthermore, it is shown in \cite[Section 3.2]{JS12} that we can parameterize $\Sigma_{KT}\cap \mathcal{U}$ by some continuous function of the $x$-variable. Hence we can prove \eqref{ktm} by the argument in Part 1. Finally, by the same result in Part 2 and $M\subseteq\Sigma_{FJ}$, the desired result \eqref{fjm} follows.     	
\end{proof}

\subsection{Verification of  PEB over FJ-condition in Case II} 

\begin{lemma}\label{lem6.3}
	Let  \eqref{BLPP} be simple at $(\bar{x},\bar{y})\in M$ in Case II. Then  PEB over FJ-condition holds at $(\bar{x},\bar{y})$ with $L>0$.
\end{lemma}

\begin{proof}
	
	In Case II, by definition, $S(\bar{x})=\{\bar{y}_1,\bar{y}_2 \}$ and there exist $C^2$-mappings $y_1(x), y_2(x)$ being unique local minimizers for $L(x)$ in a neighborhood of $\bar{x}$ with $y_1(\bar{x})=\bar{y}_1$, $y_2(\bar{x})=\bar{y}_2$ according to Type 1.  
	
	Without loss of generality we assume that $\bar{y}=\bar{y}_1$. Now we verify PEB over FJ-condition at $(\bar{x},\bar{y})$. 
	By definition, for all $v$ in $\mathcal{V}$, a small neighbourhood of $0$, 
	$
	K_{FJ}(v)\subseteq \Sigma_{FJ}\subseteq \Sigma_{gc}.
	$
	Since $(\bar x,\bar y)=(\bar x,\bar y_1)$ is of Type 1, by taking a small neighborhood $\mathcal{U}$ of $(\bar{x},\bar{y})$, points $(x,y)\in\mathcal{U}\cap  K_{FJ}(v)$ implies that  $y=y_1(x)$. Next, by \eqref{alpha}, in a neighborhood of $\bar{x}$,
	\begin{equation}\label{alpha2}
	\left|\frac{d}{dx}[f(x,y_2(x))-f(x,y_1(x))]\right|\geq\frac{|\alpha|}{2}>0.
	\end{equation}
	Clearly, we can without loss of generality assume that the above property holds on  $\mathcal{U}\cap  K_{FJ}(v)$ and $v\in \mathcal{V}$. We claim that for some large constant $L>0$, 
	\begin{equation}\label{peb}
	\mathrm{dist}_{K_{FJ}(0)}(x,y)\leq L|v|,\quad \forall\,(x,y)\in\mathcal{U}\cap K_{FJ}(v), \ v\in \mathcal{V},
	\end{equation}
	which means that the PEB over FJ-condition holds at $(\bar{x},\bar{y})$. Indeed, The proof will be divided into two cases. 
	
	\noindent
	{\bf Case 1}:  $f(x,y)-V(x)=0$.
	
	Clearly $f(x,y)-V(x)=0$ implies that $(x,y)\in K_{FJ}(0)$. Thus $\mathrm{dist}_{K_{FJ}(0)}(x,y)=0$ and the inequality \eqref{peb} holds with $L=0$.

	\noindent
	{\bf Case 2}: $f(x,y)-V(x)\neq0$.
	
	On one hand, since $y=y_1(x)$ for $(x,y)\in\mathcal{U}\cap K_{FJ}(v),v\in \mathcal{V}$ and $(\bar{x},y_1(\bar{x}))=(\bar{x},\bar{y}_1)\in K_{FJ}(0)$, for some positive constant $C$, we obtain 
	\begin{align}
	\mathrm{dist}_{K_{FJ}(0)}(x,y)=\mathrm{dist}_{K_{FJ}(0)}\left (x,y_1(x)\right )&\leq \|x-\bar{x}\|
	+\|y_1(x)-y_1(\bar{x})\|\nonumber \\ 
	&\leq C \|x-\bar{x}\|,\label{distance}
	\end{align}
	since $y_1(x)$ is a $C^2$ mapping. On the other hand, since $f(x,y_1(x))-V(x)\neq0$, by Lemma~\ref{lem1} and $S(\bar{x})=\{\bar{y}_1,\bar{y}_2 \}$, the argument in Part 1 of the proof of Lemma \ref{lem6.2} implies that $V(x)=f(x,y_2(x))$. Thus, for all $(x,y)\in\mathcal{U}\cap K_{FJ}(v), v\in \mathcal{V}$, 
	\begin{align*}
	f(x,y)-V(x)&=f(x,y_1(x))-f(x,y_2(x))\\
	&=[f(x,y_1(x))-f(x,y_2(x))]-[f(\bar{x},y_1(\bar{x}))-f(\bar{x},y_2(\bar{x}))]\\
	&=\frac{d}{dx}[f(x,y_1(x))-f(x,y_2(x))]\Big|_{x=\widetilde{x}}(x-\bar{x}),
	\end{align*}
	where $\widetilde{x}$ is some point between $x$ and $\bar{x}$ obtained by the mean value theorem. Hence by \eqref{alpha2}, since $f(x,y)-V(x)\geq0$, we have 
	\begin{equation} \label{distance1}
	f(x,y)-V(x)\geq \frac{|\alpha|}{2} \|x-\bar{x}\|, \  \mathrm{if}\ f(x,y)-V(x)\neq0.
	\end{equation}
	Hence combining (\ref{distance}) and (\ref{distance1}), for positive constant $L=\frac{2C}{|\alpha|}>0$, we have 
	\begin{equation*}
	\mathrm{dist}_{K_{FJ}(0)}(x,y)\leq L\left(f(x,y)-V(x)\right),
	\end{equation*}
	which proves \eqref{peb} and then completes the proof. 
\end{proof}

\section{Application to optimality conditions}\label{sec7}

When the lower level constraint $Y(x)=Y$ is independent of $x$, and the partial calmness constant $\mu=0$ or $\mu>0$ but the value function $V(x)$ is continuously differentiable at the point of interest, for problem \eqref{CPB1}, we can use the new constraint qualification  given in \cite[Theorem 5]{GY17} and  the new sharp necessary optimality condition given in \cite[Theorem 6]{GY19} to derive a necessary optimality condition. Similarly to the discussion in \cite{GY19},  this necessary optimality condition would be  sharper than the usual M-stationary condition based on the value function \cite[Theorem 4.1]{Ye11} for the combined program  (CP), provided that LICQ holds for the lower level program at $\bar y$. 

Now we derive optimality conditions for the general  bilevel programming problem \eqref{BLPP}  with no upper level constraints which are  simple at $(\bar{x},\bar{y})\in M$.
First we prove that MPCC-LICQ is satisfied for the partially penalized problem in all cases.

\begin{proposition}\label{MPCCLICQ}
	Let \eqref{BLPP} be simple at $(\bar{x},\bar{y})\in M$. 
	\begin{enumerate}[(1)]
		\item If $S(\bar{x})=\{\bar{y}\}$ and  KKT condition holds at $\bar{y}$ for $L(\bar{x})$, then for any KKT-multiplier $\bar{u}$, MPCC-LICQ is satisfied {at} $(\bar x, \bar y, \bar u)$ for program \eqref{CP1} with $\mu=0$. 
		
		\item If $S(\bar{x})=\{\bar{y}\}$ and  KKT condition does not hold at $\bar{y}$ for $L(\bar{x})$, then for any FJ-multiplier $(\bar{u}_0,\bar{u})\in\mathbb{R}\times\mathbb{R}^p$ at $\bar{y}$ for $L(\bar{x})$, where $\bar{u}_0=0$, MPCC-LICQ is satisfied at $(\bar x, \bar y, \bar{u}_0,\bar u)$  for program \eqref{CP2} with $\mu=0$.
		
		\item If $S(\bar{x})=\{\bar{y}_1,\bar{y}_2\}$, where $\bar{y}_1=\bar{y}$, then for the  unique KKT-multiplier $\bar{u}$ {of $\bar{y}$ for $L(\bar{x})$}, MPCC-LICQ is satisfied {at} $(\bar x, \bar y, \bar u)$ for program  \eqref{CP1}. 
	\end{enumerate}
\end{proposition}
\begin{proof}
	Define 
	$\mathcal{J}:=\big\{j: g_j(\bar x,\bar y)=0 \big\},\ \mathcal{K}:=\big\{k: \bar u_k=0 \big\}.$
	Then by \cite{SS01},  MPCC-LICQ holds at $(\bar x,\bar y,\bar u)$ for program \eqref{CP1} with $\mu=0$ if and only if the matrix 
	\begin{align*}
	M(\mathcal{J},\mathcal{K},\bar x,\bar y,\bar u):=\begin{pmatrix}
	(\nabla_{yx}^2 f+\nabla_{yx}^2 g^T\bar u)(\bar x,\bar y) & \nabla_x g_{\mathcal{J}} (\bar x,\bar y) 
	\vspace{2mm}\\
	(\nabla_{yy}^2 f+\nabla_{yy}^2 g^T\bar u)(\bar x,\bar y)  & \nabla_y g_{\mathcal{J}} (\bar x,\bar y) 
	\vspace{2mm}\\
	\nabla_y g_{\mathcal{K}^c}(\bar x,\bar y) ^T & 0
	\end{pmatrix}
	\end{align*}
	has full column rank. {Here $g_{\mathcal{J}}$ denotes the functions with components $g_j$, $j\in\mathcal{J}$ and $\mathcal{K}^c$ denotes the complement of $\mathcal{K}$.} For simplicity we omit the dependence on $(\bar x,\bar y)$ in the proof if it is clear from the context. The proof falls naturally into four parts.
	
	\textbf{Part 1: Case I, Type 1.}
	{
		In this case, by SC, we have $\mathcal{J}=\mathcal{K}^c$. Moreover, by LICQ,  $\nabla_y g_{\mathcal{J}}(\bar x,\bar y)$ has full column rank, and $(\nabla_{yy}^2 f+\nabla_{yy}^2 g^Tu)(\bar x,\bar y)|_{T_{\bar{y}}Y(\bar{x})}$ is nonsingular by using SOC and $T_{\bar{y}}Y(\bar{x})=\ker([\nabla_y g_{\mathcal{J}}] ^T)$. Hence 
		\begin{align*}
		\begin{pmatrix}
		\nabla_{yy}^2 f+\nabla_{yy}^2 g^T{{\bar u} } & \nabla_y g_{\mathcal{J}} 
		\vspace{2mm}\\
		[\nabla_y g_{\mathcal{K}^c}]^T & 0
		\end{pmatrix}
		=
		\begin{pmatrix}
		\nabla_{yy}^2 f+\nabla_{yy}^2 g^T{{\bar u} }& \nabla_y g_{\mathcal{J}} 
		\vspace{2mm}\\
		[\nabla_y g_{\mathcal{J}}]^T & 0
		\end{pmatrix}
		\mathrm{is\ nonsingular},
		\end{align*}
		which implies $M(\mathcal{J},\mathcal{K},\bar x,\bar y,\bar u)$ has full column rank. Since $\mathcal{J}=\mathcal{K}^c$, MPCC-LICQ reduces to  LICQ in this case. }
	
	\textbf{Part 2: Case I, Type 2.}
	In this case, by A2 and A3 in \cite{JS12}, $\mathcal{J}=\{1,\dots,p\}$, $p\geq1$ and $|\mathcal{K}^c|=p-1$. As in \cite{JS12}, after renumbering we assume $\mathcal{K}^c=\{1,\dots,p-1\}$. As in Part I, by A1, A5 in \cite{JS12}, it is easy to obtain that the matrix
	\begin{align*}
	C:=\begin{pmatrix}
	\nabla_{yy}^2 f+\nabla_{yy}^2 g^T\bar u & \nabla_y g_{\mathcal{K}^c} 
	\vspace{2mm}\\
	[\nabla_y g_{\mathcal{K}^c}]^T & 0
	\end{pmatrix}
	\mathrm{is\ nonsingular}.
	\end{align*}
	Note that 
	\begin{align*}
	\begin{pmatrix}
	\nabla_{yx}^2 f+\nabla_{yx}^2 g^T \bar u & \nabla_x g_{\mathcal{J}} 
	\vspace{2mm}\\
	\nabla_{yy}^2 f+\nabla_{yy}^2 g^T \bar u & \nabla_y g_{\mathcal{J}} 
	\vspace{2mm}\\
	[\nabla_y g_{\mathcal{K}^c}]^T & 0
	\end{pmatrix}
	=\begin{pmatrix}
	\nabla_{yx}^2 f+\nabla_{yx}^2 g^T \bar u & \nabla_x g_{\mathcal{K}^c} & \nabla_x g_p 
	\vspace{2mm}\\
	\nabla_{yy}^2 f+\nabla_{yy}^2 g^T \bar u & \nabla_y g_{\mathcal{K}^c} & \nabla_y g_p 
	\vspace{2mm}\\
	[\nabla_y g_{\mathcal{K}^c}]^T & 0 & 0
	\end{pmatrix}
	=\begin{pmatrix}
	A & B\\
	C & D
	\end{pmatrix},
	\end{align*}
	where the matrix $C$ is nonsingular. To show that $M(\mathcal{J},\mathcal{K},\bar x,\bar y,\bar u)$ has full column rank, it suffices to prove that 
	\begin{align*}
	\begin{pmatrix}
	A & B\\
	C & D
	\end{pmatrix}
	\begin{pmatrix}
	\lambda \\
	\mu 
	\end{pmatrix}
	=0\ \implies 
	\begin{pmatrix}
	\lambda \\
	\mu 
	\end{pmatrix}
	=0.
	\end{align*}
	By the direct computation, 
	\begin{align*}
	\begin{pmatrix}
	A & B\\
	C & D
	\end{pmatrix}
	\begin{pmatrix}
	\lambda \\
	\mu 
	\end{pmatrix}
	=\begin{pmatrix}
	A\lambda+B\mu \\
	C\lambda+D\mu 
	\end{pmatrix}=0
	\implies A\lambda+B\mu=0\ \mathrm{and}\ C\lambda+D\mu=0.
	\end{align*}
	Since $C$ is nonsingular, $\lambda=-C^{-1}D\mu$. Thus $(B-AC^{-1}D)\mu=0$. Hence $\mu=0$ and then $\lambda=0$ if $B-AC^{-1}D$ is nonsingular. Note that $B-AC^{-1}D\in \mathbb{R}$ and $B-AC^{-1}D\neq0$ can be derived by A6 in \cite{JS12} and hence MPCC-LICQ holds.

	\textbf{Part 3: Case I, Type 4.}
	In this case  KKT-condition does not hold. For any FJ-multiplier $(\bar{u}_0,\bar{u})\in\mathbb{R}\times\mathbb{R}^p$ at $\bar{y}\in S(\bar x)$, where $\bar{u}_0=0$, MPCC-LICQ is satisfied for program \eqref{CP2} if and only if the matrix 
	\begin{align*}
M(\mathcal{J},\mathcal{K},\bar x,\bar y,\bar u_0,\bar u):=	\begin{pmatrix}
	\bar u_0\nabla_{yx}^2 f+\nabla_{yx}^2 g^T\bar u & \nabla_x g_{\mathcal{J}} & 0 & 0
	\vspace{2mm}\\
	\bar u_0\nabla_{yy}^2 f+\nabla_{yy}^2 g^T\bar u & \nabla_y g_{\mathcal{J}} & 0 & 0
	\vspace{2mm}\\
	[\nabla_y g_{\mathcal{K}^c}]^T & 0 & 1_{\mathcal{K}^c} & 0
	\vspace{2mm}\\
	[\nabla_y f]^T & 0 & 0 & 1
	\end{pmatrix}
	\end{align*}
	has full column rank. Hence MPCC-LICQ holds if and only if the matrix 
	\begin{align*}
	\begin{pmatrix}
	u_0\nabla_{yx}^2 f+\nabla_{yx}^2 g^T \bar u & \nabla_x g_{\mathcal{J}} & 0 
	\vspace{2mm}\\
	u_0\nabla_{yy}^2 f+\nabla_{yy}^2 g^T \bar u & \nabla_y g_{\mathcal{J}} & 0 
	\vspace{2mm}\\
	[\nabla_y g_{\mathcal{K}^c}]^T & 0 & 1_{\mathcal{K}^c} 
	\end{pmatrix}
	\end{align*}
	has full column rank. 
	
	In this case, by C1 and C4 in \cite{JS12}, $\mathcal{K}^c=\mathcal{J}=\{1,\dots,p\}$, $p\geq1$. Let us define the following map {{$\Phi:\mathbb{R}\times\mathbb{R}^m\times\mathbb{R}\times
			\mathbb{R}^{p}\rightarrow\mathbb{R}^m \times\mathbb{R}^p\times\mathbb{R}$}}:
	{
		\begin{align*}
		\Phi(x,y,u_0,u)=
		\begin{pmatrix}
		u_0\nabla_y f(x,y)+u\nabla_y g(x,y) \vspace{2mm}\\
		{{g(x,y)} }\vspace{2mm}\\
		\sum_{j=0}^p u_j-1 
		\end{pmatrix}.
		\end{align*}}
	Note that the map $\Phi$ has the same property of the map $\Psi$ defined in \cite{JS12}. The only difference is that: the authors in \cite{JS12} normalize $\bar{u}_j$'s ($\bar{\mu}_j$'s in \cite{JS12}) by setting $\bar{\mu}_p=1$; we normalize $\bar{u}_j$'s by setting $\sum_{j=0}^p u_j=1$. Since $\bar{u}_j\neq0$ by C4 in \cite{JS12}, they are equivalent. Due to C5, $D_{x,y,u}\Psi(\bar{x},\bar{y},0, \bar{u})$ and then $D_{x,y,u}\Phi(\bar{x},\bar{y},0, \bar{u})$ is nonsingular. Note that 
	\begin{align*}
	D_{x,y,u}\Phi(\bar{x},\bar{y},0,\bar{u})=
	\begin{pmatrix}
	\bar u_0\nabla_{yx}^2 f+\nabla_{yx}^2 g^T \bar u & \nabla_x g_{\mathcal{J}} & 0 
	\vspace{2mm}\\
	\bar u_0\nabla_{yy}^2 f+\nabla_{yy}^2 g^T \bar u & \nabla_y g_{\mathcal{J}} & 0 
	\vspace{2mm}\\
	[\nabla_y g_{\mathcal{K}^c}]^T & 0 & 1_{\mathcal{K}^c} 
	\end{pmatrix}.
	\end{align*}
	Hence MPCC-LICQ holds. Moreover, since $\mathcal{J}=\mathcal{K}^c$, it reduces to LICQ in this case. 
	
	\textbf{Part 4: Case I, Type 5.}
	In this case, by D1, $\mathcal{J}=\{1,\dots,p\}$, $p=m+1\geq2$. Although MFCQ may does not hold in Type 5, KKT condition holds. We show that MPCC-LICQ holds for all Lagrange multipliers $\bar u$. To this end, we write the matrix 
	\begin{align*}
	\begin{pmatrix}
	\nabla_{yx}^2 f+\nabla_{yx}^2 g^T\bar u & \nabla_x g_{\mathcal{J}} 
	\vspace{2mm}\\
	\nabla_{yy}^2 f+\nabla_{yy}^2 g^T\bar u & \nabla_y g_{\mathcal{J}} 
	\vspace{2mm}\\
	[\nabla_y g_{\mathcal{K}^c}]^T & 0
	\end{pmatrix}
	=\begin{pmatrix}
	A & B\\
	C & 0
	\end{pmatrix},
	\end{align*}
	where $B:=\nabla g_{\mathcal{J}} $ and $C:=[\nabla_y g_{\mathcal{K}^c}]^T
	$.
	To prove the matrix $M(\mathcal{J},\mathcal{K},\bar x,\bar y,\bar u)$ has full column rank, it suffices to show that both $B$ and $C$ have full column rank. By D1 and D2 in \cite{JS12}, the matrix $B$ has full column rank and then it is nonsingular. This implies that $B$ is also full row rank. In particular, $\begin{pmatrix}\nabla_y g_{\mathcal{J}}
	\end{pmatrix}$ has full row rank. Since $\mathcal{K}^c\subseteq \mathcal{J}$, the matrix $C$ has full column rank. Therefore, $M(\mathcal{J},\mathcal{K},\bar x,\bar y,\bar u)$ has full column rank and then MPCC-LICQ holds in this case. 
\end{proof}

Since we have shown that the partial calmness holds and MPCC-LICQ is satisfied for the partially penalized problem of the combined program in the simple case,  we can now write down the optimality conditions for the combined program and hence obtain the optimality conditions  in the following theorem.

\begin{theorem}[\bf Optimality conditions for the simple BLPPs]\label{condition1}
	
		Let problem \eqref{BLPP} be simple at its local minimizer $(\bar{x},\bar{y})\in M$ and denote by $J_0:=J_0(\bar{x},\bar{y})$. Then the following necessary optimality conditions hold.
	\begin{enumerate}[(1)]
		\item[{\rm (1)}] Suppose  $S(\bar{x})=\{\bar{y}\}$ and $(\bar{x},\bar{y})$ is of Type 1.  Let $\bar{u}$ be the unique KKT-multiplier for $L(\bar{x})$ at $\bar y$. Then $(\bar{x},\bar{y},\bar{u})$ is a local solution of the following standard nonlinear programming problem (NLP):
		\begin{equation}\label{(6.1)}
		\begin{aligned}
		\min_{x,y,u}\ & F(x,y) \\
		\mathrm{s.t.\ \,}
		&\nabla_y f(x,y)+\sum_{j\in J_0} u_j\nabla_y g_j(x,y)=0, \ 
		g_j(x,y)=0,\ \forall\,j\in J_0.
		\end{aligned}
		\end{equation} 		
		And  KKT condition holds for problem (\ref{(6.1)}) at $(\bar{x},\bar{y},\bar{u})$, i.e., there exist $w\in\mathbb{R}^m$, $\xi\in\mathbb{R}^p$ such that 
		\begin{eqnarray}
		&&0=\nabla_x F(\bar{x},\bar{y}) -\big[\nabla_{xy}^2 (f+\bar{u}^Tg)(\bar{x},\bar{y})\big]w+\nabla_x (\xi^Tg)(\bar{x},\bar{y}),\label{NPKKT1} \\
		&&0=\nabla_y F(\bar{x},\bar{y}) -\big[\nabla_{yy}^2 (f+\bar{u}^Tg)(\bar{x},\bar{y})\big]w
		+\nabla_y (\xi^Tg)(\bar{x},\bar{y}), \label{NPKKT2}\\
		&&\xi_i=0,\ \forall\,i\notin J_0, \label{NPKKT3}\\
		&&  \nabla_y g_j(\bar{x},\bar{y})^Tw=0,\ \forall\,j\in J_0.\nonumber 
		\end{eqnarray}
	
		\item[{\rm (2)}] Suppose $S(\bar{x})=\{\bar{y}\}$ and $(\bar{x},\bar{y})$ is of Type 2. Let $\bar{u}$ be the unique KKT-multiplier  for $L(\bar{x})$ at $\bar y$.  Then there exists a unique $q\in J_0$ such that $\bar{u}_{q}=0$, $\bar{u}_j>0$ for $j\in J_0\setminus\{q\}$, and  $(\bar{x},\bar{y},\bar{u})$ is a local minimizer of the following MPCC:
		\begin{equation}\label{pcase2}
		\begin{aligned}
		\min_{x,y,u}\ & F(x,y) \\
		\mathrm{s.t.\ \,}
		&\nabla_y f(x,y)+\sum_{j\in J_0} u_j\nabla_y g_j(x,y)=0, \ 
		g_j(x,y)=0,\ \forall\,j\in J_0\setminus \{q\},\\
		&g_{q}(x,y) \leq 0,\  u_{q} \geq 0,\ 
		u_{q} g_{q}(x,y) =0.
		\end{aligned}
		\end{equation} 	
		And the Strong (S-) stationary condition holds for problem (\ref{pcase2}) at $(\bar{x},\bar{y},\bar{u})$, i.e., there exists $(w,\xi)\in\mathbb{R}^m\times\mathbb{R}^p$ such that  (\ref{NPKKT1})-(\ref{NPKKT3}) and the following condition holds:
		\begin{equation}\label{case2}
\nabla_y g_j(\bar{x},\bar{y})^Tw=0,\ \forall\,j\in J_0{(\bar{x},\bar{y})}\setminus\{q\},\ \nabla_y g_{q}(\bar{x},\bar{y})^T w\leq0,\ \xi_{q}\geq0.
		\end{equation} 
		
		\item[{\rm (3)}] Suppose  $S(\bar{x})=\{\bar{y}\}$ and $(\bar{x},\bar{y})$ is of Type 4.
		Then $M_{\rm FJ}(\bar x,\bar y)=\{(\bar u_0, \bar u)\}$ with $\bar u_0=0$  and  $\bar u_j>0,\ \forall\,j \in J_0$. Moreover,  $(\bar{x},\bar{y},0,\bar{u})$ is a local minimizer of the following standard NLP:
		\begin{equation}\label{pcase3}
		\begin{aligned}
		\min_{x,y,u_0,u}\ & F(x,y) \\
		\mathrm{s.t.\ \,}
		& u_0\nabla_y f(x,y)+\sum_{j\in J_0} u_j\nabla_y g_j(x,y)=0, \\
		&g_j(x,y)=0,\ \forall\,j\in J_0,\ 
		u_0\geq 0,\ u_0+\sum_{j\in J_0} u_j=1.
		\end{aligned}
		\end{equation} 
And  KKT condition holds for problem (\ref{pcase3}) at $(\bar{x},\bar{y},0, \bar{u})$, i.e., 		
		there exist  $w\in\mathbb{R}^m$, $\xi\in\mathbb{R}^p$ such that   (\ref{NPKKT1})-(\ref{NPKKT3}) without the terms involving $f$ and the following condition holds:
		\begin{equation}\label{case3}
		\nabla_y f(\bar{x},\bar{y})^T w\leq0,\ \nabla_y g_j(\bar{x},\bar{y})^T w=0,\ \forall\,j\in J_0(\bar{x},\bar{y}).
		\end{equation}
		
		\item[{\rm (4)}] Suppose $S(\bar{x})=\{\bar{y}\}$ and $(\bar{x},\bar{y})$ is of Type 5-1. Then there exists a unique KKT-multiplier $\bar{u}$ and an index $q\in J_0$ with $\bar{u}_q=0$, $\bar{u}_j>0,\ \forall\,j\in J_0\setminus\{q\}$,  $\bar{\mu}_j>0$, $\forall\,j\in J_0$ such that the set of KKT-multipliers have the representation {
			\begin{equation}\label{kkttype51}
			M_{\rm KKT}(\bar{x},\bar{y})=\bar{u}+\mathbb{R}_+ \bar{\mu}.
			\end{equation}
	 Point $(\bar{x},\bar{y},\bar{u})$ is a local solution of MPCC \eqref{pcase2} and S-stationary condition holds for problem (\ref{pcase2}) at $(\bar{x},\bar{y},\bar{u})$ with $w=0$. Furthermore, $(\bar{x},\bar{y})$ is a local solution of the following standard NLP:
		\begin{equation}\label{pcase4+}
		\begin{aligned}
		\min_{x,y}\ & F(x,y) \\
		\mathrm{s.t.\ \,}
		&g_j(x,y)=0,\ \forall\,j\in J_0{(\bar{x},\bar{y})}\setminus\{q\},\ 
		g_{q}(x,y)\leq0,
		\end{aligned}
		\end{equation} and KKT condition for the above NLP  holds i.e.,
			there exists $\xi\in \mathbb{R}^p$ such that
			\begin{equation}\label{type51xi}
			\begin{aligned}
			&0=\nabla F(\bar{x},\bar{y}) +\nabla (\xi^Tg)(\bar{x},\bar{y}),\\
			&\xi_i=0,\ \forall\,i\notin J_0(\bar{x},\bar{y}),\ 
			\xi_{q}\geq 0.
			\end{aligned} 
			\end{equation} 
		}

	\item[{\rm (5)}] Suppose $S(\bar{x})=\{\bar{y}\}$ and $(\bar{x},\bar{y})$ is of Type 5-2. Then there exists a unique pair of $q,r\in J_0$, $q\neq r$ such that there are  $\bar{u}^1,\bar{u}^2\in M_{\rm KKT}(\bar{x},\bar{y})$ satisfying 
		\begin{equation*}\label{type52kkt}
		\bar{u}^1_{q}=0,\ \bar{u}_j^1>0,\ \forall\,j\in J_0\setminus\{q\};\quad 
		\bar{u}^2_{r}=0,\ \bar{u}_j^2>0,\ \forall\,j\in J_0\setminus\{r\},
		\end{equation*}
		and $M_{\rm KKT}(\bar{x},\bar{y})$ is the  line segment jointing $\bar{u}^1$ and $\bar{u}^2$.
		Moreover, for any $\bar{u}\in M_{\rm KKT}(\bar{x},\bar{y})$, $(\bar{x},\bar{y},\bar{u})$ is a local solution of the following MPCC:	
		{
		\begin{equation}\label{pcase5+}
		\begin{aligned}
		\min_{x,y,u}\ & F(x,y) \\
		\mathrm{s.t.\ \,}
		&\nabla_y f(x,y)+\sum_{j\in J_0} u_j\nabla_y g_j(x,y)=0, \ 
		g_j(x,y)=0,\ \forall\,j\in J_0\setminus \{q,r\},\\
		& g_{k}(x,y) \leq 0,\  u_{k} \geq 0,\ u_{k} g_{k}(x,y) =0,\ \forall\,k\in\{q,r\}.
		\end{aligned}
		\end{equation}	
	And S-stationary condition holds for problem (\ref{pcase5+}) at $(\bar{x},\bar{y},\bar{u})$. Moreover, S-stationary condition reduces to the existence of  $\xi\in\mathbb{R}^{p}$ such that 
	\begin{equation*}
			\begin{aligned}
			&0=\nabla F(\bar{x},\bar{y}) +\nabla (\xi^Tg)(\bar{x},\bar{y}),\\
			&\xi_i=0,\ \forall\,i\notin J_0(\bar{x},\bar{y}),\ 
			\xi_{q}\geq 0 \mbox{ if } \bar{u}=\bar{u}^1,\  \xi_{r}\geq 0 \mbox{ if } \bar{u}=\bar{u}^2 .
			\end{aligned} 
			\end{equation*} 

}

		
		\item[{\rm (6)}] Case II: $S(\bar{x})=\{\bar{y}_1,\bar{y}_2\}$. Let $y_1(x),y_2(x)$ be the unique local minimizers for $L({x})$ in a neighborhood of $\bar x$ with $y_k(\bar x)=\bar{y}_k$, and $\bar{u}^k$ be the unique KKT-multiplier for $L(\bar{x})$ at $\bar{y}_k$, $k=1,2$. Then the value function can be written as 
		$
		V(x)=\min\big\{f(x,y_1(x)),f(x,y_2(x))\big\}
		$
		for all $x$ lying in a neighbourhood of $\bar x$. Let $\bar{u}:=\bar{u}^1$ be the unique KKT-multiplier for $L(\bar{x})$ at $\bar{y}=\bar{y}_1$. Then $\bar{u}_j>0$ for all $j\in J_0(\bar{x},\bar{y})$ and $(\bar{x},\bar{y},\bar{u})$ is a local solution of the following nonsmooth NLP for some $\mu > 0$:
		\begin{equation}\label{pcase 6}
	\begin{aligned}
	\min_{x,y,u}\ & F(x,y)+\mu\big(f(x,y)-V(x)\big) \\
	\mathrm{s.t.\ \,}
	& \nabla_y f(x,y)+u\nabla_y g(x,y)=0, \ 
	 g_j(x,y)=0,\ \forall\,j\in J_0{(\bar{x},\bar{y})}.
	\end{aligned}	
	\end{equation} 	
	And the nonsmooth KKT condition holds for problem (\ref{pcase 6}) at	 $(\bar{x},\bar{y},\bar{u})$, that is, there exist $\lambda^k\geq0$ with $\sum_{k=1}^2\lambda^k=1$,  $w\in\mathbb{R}^m$, $\xi\in\mathbb{R}^p$, and $\mu>0$ such that 
		\begin{equation*}\label{case4}
		\begin{aligned}
		&0=\nabla_x F(\bar{x},\bar{y}) 
		+\mu\sum_{k=1}^2\lambda^k\big(\nabla_x f(\bar{x},\bar{y})-\nabla_x f(\bar{x},\bar{y}_k)-\bar{u}^k\nabla_x g(\bar{x},\bar{y}_k)\big)\\   	
		&\qquad-\big[\nabla_{xy}^2 (f+\bar{u}^Tg)(\bar x,\bar y)\big]w
		+\nabla_x (\xi^Tg)(\bar{x},\bar{y}),\\		
		&0=\nabla_y F(\bar{x},\bar{y})+\mu\nabla_y f(\bar{x},\bar{y})-\big[\nabla_{yy}^2 (f+\bar{u}^Tg)(\bar x,\bar y)\big]w
		+\nabla_y (\xi^Tg)(\bar{x},\bar{y}),\\
		&\xi_i=0,\ \forall\,i\notin J_0(\bar{x},\bar{y}),\ \nabla_y g_j(\bar{x},\bar{y})^Tw=0,\ \forall\,j\in J_0(\bar{x},\bar{y}).
		\end{aligned}
		\end{equation*}  
	\end{enumerate}
\end{theorem}	
{
\begin{proof}
	For Case I, by Lemma~\ref{lem6.2}, $\Sigma_{FJ}=M$ in a neighbourhood of $(\bar{x},\bar{y})$. Hence problem \eqref{CP4} is partially calm at $(\bar x,\bar y)$ with $\mu=0$. By Theorem~\ref{thmequiv}, problem \eqref{CPFJ} is partially calm at $(\bar x,\bar y,\bar u_0, \bar u)$ for each FJ-multiplier $(\bar u_0, \bar u)\in M_{\rm FJ}(\bar{x},\bar{y})$. When $\bar{u}_0\neq0$, this means that problem \eqref{CP} is partially calm at $(\bar{x},\bar{y},\bar{u})$ for each KKT-multiplier $\bar u\in M_{\rm KKT}(\bar{x},\bar{y})$ with the same calmness modulus $\mu=0$. For Case II, by Lemma~\ref{lem6.3} and the Lipschitz continuity of the objective function $F(x,y)$, problem \eqref{CP} is partially calm at $(\bar{x},\bar{y},\bar{u})$ with a  calmness modulus $\mu>0$.
	
	By dropping inactive inequalities, it is easy to see that locally, the partially penalized problem can be reduced to corresponding optimization problem in each case.
%

	Since MPCC-LICQ (which reduces to a standard constraint qualification if the problem is a NLP) is proved in Proposition~\ref{MPCCLICQ}, S-stationary condition  hold for each optimization problem and so no extra constraint qualifications are required.

	Next, we consider the structure of FJ and KKT-multipliers. Except Type 5, the FJ and KKT-multipliers are unique. In Case I, Type 5-1, since $\big|J_0(\bar{x},\bar{y})\big|=m+1$ by (D1) in \cite{JS12} and any $m$ elements of $\big\{\nabla_{y} g_j(\bar{x},\bar{y}),\ j\in J_0\big\}$ is linearly independent by (D3) in \cite{JS12}, CRCQ holds. But MFCQ does not hold since all $\mu_j$, $j\in J_0$ in \cite{JS12} have the same sign. Since $\mu_j$ in \cite{JS12} is unique up to a common multiple, there exists a unique KKT-multiplier $\bar{u}$ and an index $q\in J_0$ such that \eqref{kkttype51} holds. On the other hand, the optimality condition \eqref{type51xi} follows from (\ref{NPKKT1})-(\ref{NPKKT3}) and \eqref{case2} since 
	\begin{equation}
	\nabla_y g_j(\bar{x},\bar{y})^Tw=0,\ \forall\,j\in J_0{(\bar{x},\bar{y})}\setminus\{q\}
	\implies w=0,
	\end{equation}
	because $\big|J_0(\bar{x},\bar{y})\setminus\{q\}\big|=m$ and $\big\{\nabla_{y} g_j(\bar{x},\bar{y}),\ j\in J_0\setminus\{q\}\big\}$ is linearly independent. In the same manner we can see that $w=0$ in Case I, Type 5-2. We could say more on Case I, Type 5-1. In this case, since $\bar{u}_j>0,\ \forall\,j\in J_0\setminus\{q\}$ for all KKT-multiplier $\bar{u}\in M_{\rm KKT}(\bar{x},\bar{y})$ and $\big\{\nabla_{y}g_j(\bar{x},\bar{y}),\ j\in J_0\setminus\{q\}\big\}$ is linearly independent in $\mathbb{R}^m$ with $\big|J_0(\bar{x},\bar{y})\setminus\{q\}\big|=m$, by the implicit function theorem, the curve $y^q(x)$ defined by  $y^q(\bar{x})=\bar{y}$ and $g_j(x,y^q(x))=0,\ \forall\,j\in J_0\setminus\{q\}$ satisfies $(x,y^q(x))\in M$ in a neighbourhood of $\bar{x}$. Hence, $(\bar{x},\bar{y})$ is actually a local solution of  NLP (\ref{pcase4+}).
		In Case I, Type 5-2, both CRCQ and MFCQ hold since any $m$ elements of $\big\{\nabla_{y} g_j(\bar{x},\bar{y}),\ j\in J_0\big\}$ is linearly independent by (D1-D3) in \cite{JS12} and $\mu_j$, $j\in J_0$ have different signs. Hence, by a similar argument of the proof of   \cite[Lemma 5.2]{BT05}, there exists a unique pair of $q,r\in J_0$, $q\neq r$ such that there are  $\bar{u}^1,\bar{u}^2\in M_{\rm KKT}(\bar{x},\bar{y})$ satisfying \eqref{type52kkt} and $M_{\rm KKT}(\bar{x},\bar{y})$ is the line segment jointing $\bar{u}^1$ and $\bar{u}^2$.
\end{proof}

}

By using the implicit function theorem, we can show that Theorem \ref{condition1} and  the necessary optimality conditions derived in \cite[Theorem 4.3]{JS12}  are equivalent in the following sense.	Due to the page limit we omit the proof.
\begin{theorem}[\bf Equivalence of optimality conditions derived in Theorem 4.3 of \cite{JS12} and Theorem \ref{condition1}]\label{eqivthm}
	
{\bf Case I, Type 1:} The optimality condition in Theorem~\ref{condition1} (1) holds if and only if there exist $C^2$-mappings $y(x)$, $u_j(x)$, $j\in \mathcal{J}:=J_0(\bar{x},\bar{y})$, satisfying
\begin{equation}\label{initialpt}
y(\bar{x})=\bar{y},\quad u(\bar{x})=\bar{u},
\end{equation}
and the following equations in a neighborhood of $\bar x$: 
	\begin{equation}\label{c1t1}
	\left\{
	\begin{aligned}
	&\ \nabla_y f(x,y(x))+\sum_{j\in \mathcal{J}}u_j(x)\nabla_y g_j(x,y(x))=0, \\
	&\ \quad g_j(x,y(x))=0,\ \forall\,j\in \mathcal{J},
	\end{aligned}
	\right.
	\end{equation} such that $\nabla_x F(\bar{x},\bar{y})+\nabla_y F(\bar{x},\bar{y})\cdot D_x y(\bar{x})=0$. 

{\bf Case I, Type 2:} The optimality condition in Theorem~\ref{condition1} (2) holds if and only if there exist 
$C^2$-mappings $\widetilde{y}(x)$, $\widetilde{u}_j(x)$, $j\in J_0(\bar{x},\bar{y})\setminus\{q\}$ and $\widehat{y}(x)$, $\widehat{u}_j(x)$, $j\in J_0(\bar{x},\bar{y})$ satisfying \eqref{initialpt} and  \eqref{c1t1} with $\mathcal{J}=J_0(\bar{x},\bar{y})\setminus\{q\}$ and $\mathcal{J}=J_0(\bar{x},\bar{y})$, respectively, such that
\begin{align*}
	\left\{
	\begin{array}{ll}
	\nabla_x F(\bar{x},\bar{y})+\nabla_y F(\bar{x},\bar{y})\cdot D_x \widetilde{y}(\bar{x})\leq0  ,\\
	\nabla_x F(\bar{x},\bar{y})+\nabla_y F(\bar{x},\bar{y})\cdot D_x \widehat{y}(\bar{x})\geq0  ,
	\end{array}
	\right.
	\mathrm{if}\ \frac{d g_{q}(x,\widetilde{y}(x))}{dx}\Big|_{x=\bar{x}}> 0,\\
	\left\{
	\begin{array}{ll}
	\nabla_x F(\bar{x},\bar{y})+\nabla_y F(\bar{x},\bar{y})\cdot D_x \widetilde{y}(\bar{x})\geq0  ,\\
	\nabla_x F(\bar{x},\bar{y})+\nabla_y F(\bar{x},\bar{y})\cdot D_x \widehat{y}(\bar{x})\leq0, 
	\end{array}
	\right.
	\mathrm{if}\ \frac{d g_{q}(x,\widetilde{y}(x))}{dx}\Big|_{x=\bar{x}}< 0.
	\end{align*}
	
	{\bf Case I, Type 4:} The optimality condition in Theorem~\ref{condition1} (3) holds if and only if there exist $C^2$-mappings $x(u_0),\,y(u_0),\, u_j(u_0)$, $j\in \mathcal{J}:=J_0(\bar{x},\bar{y})$ satisfying $x(0)=\bar{x},\,y(0)=\bar{y},\,u(0)=\bar{u}$ and the following system in an open neighborhood of $\bar{u}_0$: 
	\begin{equation*}
	\left\{
	\begin{aligned}
	&\ u_0\nabla_y f(x(u_0),y(u_0))+\sum_{j\in \mathcal{J}}u_j(u_0)\nabla_y g_j(x(u_0),y(u_0))=0, \\
	&\ u_0+\sum_{j\in \mathcal{J}}u_j(u_0)=1,\\
	&\ g_j(x(u_0),y(u_0))=0,\ \forall\,j\in \mathcal{J},
	\end{aligned}
	\right.
	\end{equation*} such that 
	$
	\nabla_x F(\bar{x},\bar{y})D_{u_0} x(0)+\nabla_y F(\bar{x},\bar{y})\cdot D_{u_0} y(0)
\geq 0.
	$
	
		{\bf Case I, Type 5-1:} The optimality condition in Theorem~\ref{condition1} (4) holds if and only if there exist $C^2$-mappings $y^q(x)$, $u_j^q(x)$, $j\in J_0(\bar{x},\bar{y})\setminus\{q\}$ satisfying \eqref{initialpt} and \eqref{c1t1} with $\mathcal{J}=J_0(\bar{x},\bar{y})\setminus\{q\}$ such that
		\begin{align*}
	\left\{
	\begin{array}{ll}
	\nabla_x F(\bar{x},\bar{y})+\nabla_y F(\bar{x},\bar{y})\cdot D_x y^q(\bar{x})\leq0  &\mathrm{if}\ \frac{d g_q(x,y^q(x))}{dx}\big|_{x=\bar{x}}> 0,\\
	\nabla_x F(\bar{x},\bar{y})+\nabla_y F(\bar{x},\bar{y})\cdot D_x y^q(\bar{x})\geq0  &\mathrm{if}\ \frac{d g_q(x,y^q(x))}{dx}\big|_{x=\bar{x}}< 0.
	\end{array}
	\right.
	\end{align*}
	
		{\bf Case I, Type 5-2:} The optimality condition in Theorem~\ref{condition1} (5) holds if and only if there exist $C^2$-mappings curves $y^q(x), u_j^q(x)$, $j\in J_0(\bar{x},\bar{y})\setminus\{q\}$, $y^r(x), u_j^r(x)$, $j\in J_0(\bar{x},\bar{y})\setminus\{r\}$ satisfying \eqref{initialpt} with $\bar{u}=\bar{u}^1, \bar{u}^2$ and \eqref{c1t1} with $\mathcal{J}=J_0(\bar{x},\bar{y})\setminus\{q\}$ and $\mathcal{J}=J_0(\bar{x},\bar{y})\setminus\{r\}$, respectively, such that
	\begin{align}\label{c1t52n}
	\left\{
	\begin{array}{ll}
	\nabla_x F(\bar{x},\bar{y})+\nabla_y F(\bar{x},\bar{y})\cdot D_x {y}^q(\bar{x})\leq0  ,\\
	\nabla_x F(\bar{x},\bar{y})+\nabla_y F(\bar{x},\bar{y})\cdot D_x {y}^r(\bar{x})\geq0  ,
	\end{array}
	\right.
	\mathrm{if\ }\gamma_q>0\ 
	\left(\mathrm{i.e.\ }\gamma_r<0\right),\\
	\left\{
	\begin{array}{ll}
	\nabla_x F(\bar{x},\bar{y})+\nabla_y F(\bar{x},\bar{y})\cdot D_x {y}^q(\bar{x})\geq0  ,\\
	\nabla_x F(\bar{x},\bar{y})+\nabla_y F(\bar{x},\bar{y})\cdot D_x {y}^r(\bar{x})\leq0  ,
	\end{array}
	\right.
	\mathrm{if\ }\gamma_q<0\ 
	\left(\mathrm{i.e.\ }\gamma_r>0\right),\label{c1t52n+}
	\end{align}
	where $\gamma_q:=\mathrm{sign}\left(\frac{d g_q(x,y^q(x))}{dx}\Big|_{x=\bar{x}}\right)$ and $\gamma_r:=\mathrm{sign}\left(\frac{d g_r(x,y^r(x))}{dx}\Big|_{x=\bar{x}}\right)$.
	
		{\bf Case II:} The optimality condition in Theorem~\ref{condition1} (6) holds if and only if there exist $C^2$-mappings $y_k(x)$, $u_j^k(x)$, $j\in J_0(\bar{x},\bar{y}_k)$, $k=1,2$, satisfying \eqref{initialpt} with $\bar{y}=\bar{y}_k, \bar{u}=\bar{u}^k$ and  \eqref{c1t1} with $\mathcal{J}=J_0(\bar{x},\bar{y}_k)$ such that 
		\begin{align*}
	\left\{
	\begin{array}{ll}
	\nabla_x F(\bar{x},\bar{y})+\nabla_y F(\bar{x},\bar{y})\cdot D_x y_1(\bar{x})\geq0  &\mathrm{if}\ \frac{d\left[f(x,y_2(x))-f(x,y_1(x))\right]}{dx}\big|_{x=\bar{x}}> 0,\\
	\nabla_x F(\bar{x},\bar{y})+\nabla_y F(\bar{x},\bar{y})\cdot D_x y_1(\bar{x})\leq0  &\mathrm{if}\ \frac{d\left[f(x,y_2(x))-f(x,y_1(x))\right]}{dx}\big|_{x=\bar{x}}< 0.
	\end{array}
	\right.
	\end{align*}
\end{theorem}

\begin{remark}
	This result is related to an interesting and involved issue in  \cite[Remark 4.2]{JS12}: the comparison of optimality conditions derived in the literature \cite[Theorem 4.3]{JS12} and from the literature \cite{DD12,YZ95,YZ10} by using KKT, value function, and the combined approaches. From our point of view, the underlying reason for the equivalence in Theorem~\ref{eqivthm} is that the structure of FJ-multipliers can distinguish the local structure of bilevel problems. The generic structure of FJ-multipliers locally around global minimizers are given in the following Figure~\ref{FJtypes}. It is worth noting that Figure~\ref{FJtypes} distinguishes all of the cases. 
	

	\begin{figure}[h!]
		\centering
		
		\includegraphics[scale=0.4]{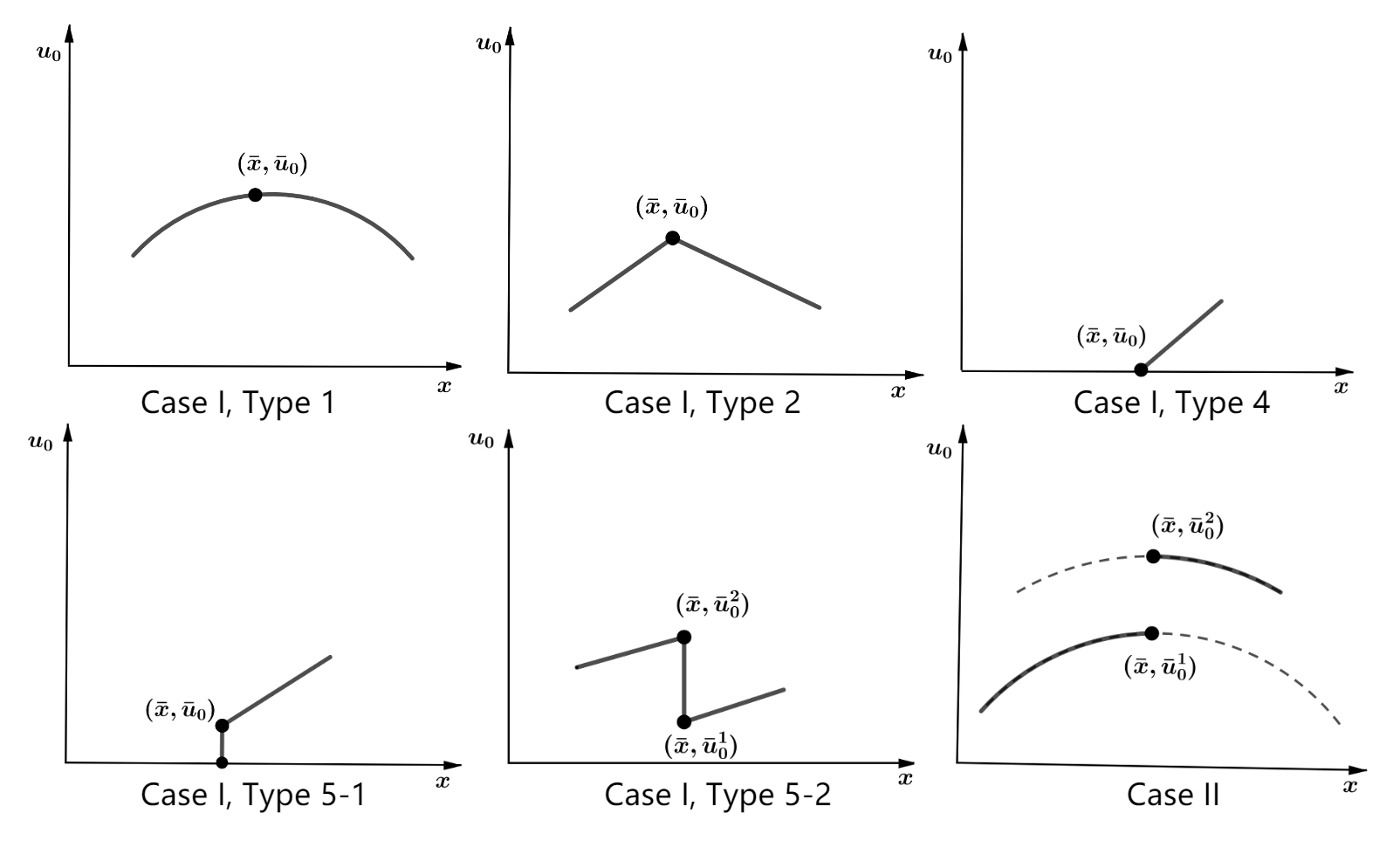}
		\caption{Generic structure of FJ-multipliers locally around global minimizers.}
		\label{FJtypes}
	\end{figure}
\end{remark}

{The following  example is a typical example of  Case I, type 4. Let us demonstrate our results by deriving optimality conditions for this example.}
\begin{example}[\text{\cite[Example 2.2]{JS12}}]\label{eg22}
	\begin{align*}
	\min_{x,y}\ & x+\sum_{j=1}^m y_j\quad 
	\mathrm{s.t.}\quad y\in S(x):=\mathrm{arg}\min_y\Big\{-y_1 : \sum_{j=1}^{m}y_j^2-x\leq0\Big\}.
	\end{align*}
	It is easy to see that $S(x)=\{(\sqrt{x},0_{m-1})\}$ if $x\geq0$ and $S(x)=\emptyset$ if $x<0$. The optimal solution of the bilevel problem is $(0,0_m)$. The lower level problem violates MFCQ at $0_m$ and $0_m$ is not a KKT-point.  One may consider whether we can use the value function approach. When $m=1$, it is not difficulty to check that $L(x)$ has a uniformly weak sharp minimum {over $f$} around $(0,0)$ and hence the value function approach can be used if $m=1$.  However since MFCQ fails, the value function $V( x)=\sqrt{x}$ is not Lipschitz continuous at $\bar x=0$ and it is difficult to obtain a necessary optimality condition for problem \eqref{VP}. Moreover for $m\geq2$, the partial calmness for (VP) (i.e., partial calmness over $f$) does not hold at $(0,0_m)$. It is worth noting that the violation of the partial calmness for (VP) cannot be rectified by a small perturbation of the data of the problem. However, it is easy to check  the only FJ multiplier for problem $L(\bar x)$ is $(\bar u_0,\bar u)=(0,1)$. The optimality condition from Theorem \ref{condition1}(2) is verified with $\xi=1$ and {$w=\frac{1}{2}(1,\dots, 1)^T$.} 
	
\end{example} 

Next we consider a special case where there is no  upper and lower level constraints. That is, we consider the unconstrained  BLPPs: 
\begin{equation}\label{SBLPP}\tag{UBLPP}
\begin{aligned}
\min_{x,y}&\  F(x,y)\quad 
\mathrm{s.t.}&\ y\in \mathrm{arg}\min_y f(x,y).
\end{aligned}
\end{equation}
It is clear that only Case I, Type 1 and Case II will happen if \eqref{SBLPP} is simple at $(\bar{x},\bar{y})\in M$. Hence the following corollary follows from Theorem \ref{condition1}.
\begin{corollary}\label{thmno}
	Let the unconstrained bilevel programming problem \eqref{SBLPP} be simple at {its local minimizer} $(\bar{x},\bar{y})\in M$,  i.e., Case I type 1 or Case II, then there exist $\mu\geq0$ and $w\in\mathbb{R}^m$ such that 
	\begin{align*}
	&0=\nabla_x F(\bar{x},\bar{y})+\mu(\nabla_xf(\bar{x},\bar{y})-\nabla_xf(\bar{x},\bar{y}_2))+\nabla_{yx}^2f(\bar{x},\bar{y})w,\\
	&0=\nabla_y F(\bar{x},\bar{y})+\nabla_{yy}^2f(\bar{x},\bar{y})w,
	\end{align*}
	where $\mu=0$ in Case I and $\bar{y}_2$ is the distinct minimizer of $f(\bar{x},\cdot)$ in Case II. 
\end{corollary}

We claim that the points in Case II are $f$-critical by the notation of Mirrlees \cite{M7599}. Indeed, according to \cite{M7599}, $(\bar{x},\bar{y})$ is $f$-critical if $\bar{y}\in S(\bar{x})$, but there is no neighborhood $\mathcal{N}$ of $(\bar{x},\bar{y})$ such that $(x',y')\in\mathcal{N}$, $\nabla_y f(x',y')=0$ imply that $y'\in S(x')$. Since $\alpha\neq0$ in \eqref{alpha}, it is easily seen that $(\bar{x},\bar{y})$ is $f$-critical if $(\bar{x},\bar{y})$ is in Case II. On the other hand, if $(\bar{x},\bar{y})$ is not $f$-critical, then there is a neighborhood $\mathcal{N}$ of $(\bar{x},\bar{y})$ such that $(x',y')\in\mathcal{N}$, $\nabla_y f(x',y')=0$ imply that $y'\in S(x')$. Thus $\Sigma_{gc}\cap\mathcal{U}=M\cap\mathcal{U}$ and then {$(\bar{x},\bar{y})$} is in Case I. This implies that $f$-critical points--generically--are in Case II under our problem setting by Theorem~\ref{newthm41}. Clearly, this verifies Mirrlees' claim: ``Normally, a $V$-critical point $(a',z')$ is such that there are two distinct maxima, $z'$ and $z''$, of $V(a,\cdot)$, there being no others" on page 16 of \cite{M7599}. Here $V$ is the lower level objective function and $a,z$ are the upper and lower level variables, respectively. Moreover, under this claim, Mirrlees gave a form of KKT condition for \eqref{SBLPP} in (52) on page 17 of \cite{M7599} which coincides with the necessary optimality condition in  \cite[Theorem 2.1]{YZ10} with a rigorous proof. Hence the condition in  \cite[Theorem 2.1]{YZ10} (the one in Corollary \ref{thmno}) is, of course, a generic one.   

\section*{Acknowledgments} 
The authors would like to thank the associate editor, and the anonymous referees for their helpful suggestions and comments, which resulted in essential improvements of the paper.

\bibliographystyle{siamplain}

\begin{thebibliography}{10}

\bibitem{amsmath}
{\sc {American Mathematical Society}}, {\em User's guide for the
  \texttt{amsmath} package (version 2.0)}, 2002,
  \url{ftp://ftp.ams.org/pub/tex/doc/amsmath/amsldoc.pdf} (accessed
  2015-07-30).

\bibitem{AMSMSC2010}
{\sc {American Mathematical Society}}, {\em {Mathematics Subject
  Classification}}, 2010, \url{http://www.ams.org/mathscinet/msc/msc2010.html}
  (accessed 2015/03/29).

\bibitem{clawpack}
{\sc {Clawpack Development Team}}, {\em Clawpack software}, 2015,
  \url{http://www.clawpack.org} (accessed 2015/05/14).
\newblock Version 5.2.2.

\bibitem{CalcI}
{\sc P.~Dawkins}, {\em Paul's online math notes: Calculus {I} --- notes},
  \url{http://tutorial.math.lamar.edu/Classes/CalcI/MeanValueTheorem.aspx}
  (accessed 2015-07-08).

\bibitem{shortmath}
{\sc M.~Downes}, {\em Short math guide for {\LaTeX}}, 2002,
  \url{ftp://ftp.ams.org/pub/tex/doc/amsmath/short-math-guide.pdf} (accessed
  2015-07-30).

\bibitem{pgfplots}
{\sc C.~Feuers\"anger}, {\em Manual for package \texttt{PGFPLOTS}}, May 2015,
  \url{http://sourceforge.net/projects/pgfplots}.

\bibitem{Hi14}
{\sc N.~Higham}, {\em A call for better indexes}.
\newblock SIAM Blogs, Nov. 2014,
  \url{http://blogs.siam.org/a-call-for-better-indexes/} (accessed 2015-04-05).

\bibitem{KoMa14}
{\sc T.~G. Kolda and J.~R. Mayo}, {\em An adaptive shifted power method for
  computing generalized tensor eigenpairs}, SIAM Journal on Matrix Analysis and
  Applications, 35 (2014), pp.~1563--1581,
  \url{https://doi.org/10.1137/140951758}.

\bibitem{La86}
{\sc L.~Lamport}, {\em \LaTeX: A Document Preparation System}, Addison--Wesley,
  Reading, MA, 1986.

\bibitem{MiGo04}
{\sc F.~Mittlebach and M.~Goossens}, {\em The \LaTeX\ Companion},
  Addison--Wesley, 2nd~ed., 2004.

\bibitem{Ne03}
{\sc M.~E.~J. Newman}, {\em Properties of highly clustered networks}, Phys.
  Rev. E, 68 (2003), 026121 (6~pages),
  \url{https://doi.org/10.1103/PhysRevE.68.026121}.

\bibitem{PeKoPi14}
{\sc C.~Peng, T.~G. Kolda, and A.~Pinar}, {\em Accelerating community detection
  by using {K}-core subgraphs}, Mar. 2014,
  \url{https://arxiv.org/abs/1403.2226}.

\bibitem{WoZhMeSh05}
{\sc D.~E. Woessner, S.~Zhang, M.~E. Merritt, and A.~D. Sherry}, {\em Numerical
  solution of the {Bloch} equations provides insights into the optimum design
  of {PARACEST} agents for {MRI}}, Magnetic Resonance in Medicine, 53 (2005),
  pp.~790--799, \url{https://doi.org/10.1002/mrm.20408},
  \url{https://www.ncbi.nlm.nih.gov/pubmed/15799055} .

\bibitem{siam}
{\em {SIAM} style manual: For journals and books}, 2013,
  \url{https://www.siam.org/journals/pdf/stylemanual.pdf}.

\end{thebibliography}


\begin{thebibliography}{1}

\bibitem{KoMa14}
{\sc T.~G. Kolda and J.~R. Mayo}, {\em An adaptive shifted power method for
  computing generalized tensor eigenpairs}, SIAM Journal on Matrix Analysis and
  Applications, 35 (2014), pp.~1563--1581,
  \url{https://doi.org/10.1137/140951758}.

\end{thebibliography}


\begin{thebibliography}{10}
	
	\bibitem{AHO18} L. Adam, R. Henrion and J.V. Outrata, {\em On M-stationarity conditions in MPECs and the associated qualification conditions}, Math. Program., 168 (2018), pp. 229–259. 
	
	\bibitem{AS13} G. Allende and G. Still, {\em Solving bilevel programs with the KKT-approach}, Math. Program.,
	138 (2013), pp. 309-332.
	
	{
\bibitem{BT05} A. Baccari and A. Trad, {\em On the classical necessary second-order optimality conditions in the presence of equality and inequality constraints}, SIAM J. Optim., 15 (2005), pp. 394-408.	
}
	
	\bibitem{B98} J. Bard, {\em Practical Bilevel Optimization: Algorithms and Applications}, Dordrecht: Kluwer Academic Publishers, 1998.
	
	
	\bibitem{BJ05} M. Bjørndal and K. Jørnsten, {\em The deregulated electricity market viewed as a bilevel programming problem}, J. Global Optim., 33 (2005), pp. 465-475.
	
	
	\bibitem{BG09} P. Bond and A. Gomes, {\em Multitask principal-agent problems: optimal contracts, fragility, and effort misallocation}, J. Econ. Theory, 144 (2009), pp. 175-211.
	
	\bibitem{BL75} P. Br\"{o}cker and L. Lander, {\em Differentiable Germs and Catastrophes}, Cambridge University Press, Cambridge, 1975.
	
	\bibitem{Clarke}
	F.~Clarke,
	{\em Optimization and Nonsmooth Analysis},
	Wiley Interscience, New York, 1983.
	
	
	
	\bibitem{C09} J.R. Conlon, {\em Two new conditions supporting the first-order
		approach to multisignal principal-agent problems}, Econometrica, 77 (2009), pp. 249-278.
	
	\bibitem{D18} S. Dempe, {\em Foundations of Bilevel Programming}, Kluwer Academic Publishers, Dordrecht, 2002.
	%
	
	
	\bibitem{DD12} S. Dempe and J. Dutta, {\em Is bilevel programming a special case of a mathematical program with complementarity constraints?}, Math. Program., 131 (2012), pp. 37-48.
	
	\bibitem{DDM07} S. Dempe, J. Dutta and B.S. Mordukhovich, {\em New necessary optimality conditions in optimistic bilevel programming}, Optimization, 56 (2007), pp. 577-604.
	
	
	
	\bibitem{DGJ09} S. Dempe, H. G\"{u}nzel and H.Th. Jongen, {\em On reducibility in bilevel problems}, SIAM J. Optim., 20 (2009), pp. 718-727.
		{ \bibitem{DJS} D. Dominik, H.Th. Jongen {and} V. Shikhman, {\em On Intrinsic Complexity of Nash Equilibrium Problems and Bilevel Optimization}, J. Optim. Theory Appl., 159 (2013), pp. 606-634.}
	\bibitem{DKPK15} S. Dempe, V. Kalashnikov, G.A. P\'{e}rez-Vald\'{e}s and N. Kalashnykova, {\em Bilevel Programming Problems, Energy Systems}, Springer, Berlin. 2015.
	
	\bibitem{DZ13} S. Dempe and A.B. Zemkoho, {\em The bilevel programming problems: reformulations, constraint qualifications and optimality conditions}, Math. Program., 138 (2013), pp. 447-473.
	
	{\bibitem{Dempebook} {S. Dempe and A.B. Zemkoho},  {\em Bilevel optimization: advances and next challenges}.  Springer Optimization and its Applications, vol. 161, 2020.}
	

	
	\bibitem{FFSGP18} L. Franceschi, P. Frasconi, S. Salzo, R. Grazzi and M. Pontil, {\em Bilevel programming for hyperparameter optimization}, Proceedings of the 35th International Conference on Machine Learning, PMLR 80:1568-1577, 2018.
	
	\bibitem{GY17} H. Gfrerer and J.J. Ye, {\em New constraint qualifications for mathematical programs with equilibrium constraints via variational analysis}, SIAM J. Optim., 27 (2017), pp. 842–865.
	
	\bibitem{GY19} H. Gfrerer and J.J. Ye, {\em New sharp necessary optimality conditions for mathematical programs with equilibrium constraints}, Set-Valued Var. Anal, 28 (2020), pp. 395-426.
	
	\bibitem{GH} S.J. Grossman and O.D. Hart, {\em An analysis of the principal-agent problem}, Econometrica, 51 (1983), pp. 7-45.
	
	%
	
		\bibitem{HJO} R. Henrion, A. Jourani and J.V. Outrata, {\em On the calmness of a class of multifunctions}, SIAM J. Optim., 26 (2002), pp. 603-618. 
	\bibitem{HS11} R. Henrion and T.M. Surowiec, {\em On calmness conditions in convex bilevel programming}, Appl. Anal., 90 (2011), pp. 951-970. 
	

	
	
	\bibitem{HM87} B. Holmstr\"{o}m and P. Milgrom, {\em Aggregation and Linearity in the Provision of Intertemporal Incentives}, Econometrica, 55 (1987), pp. 303-328.
	
	\bibitem{HM91} B. Holmstr\"{o}m and P. Milgrom, {\em Multitask principal–agent analyses: Incentive contracts, asset ownership, and job design}, J. Law, Econ., Organ., 7 (1991), pp. 24–52.
	
	\bibitem{JJT86} H.Th. Jongen, P. Jonker and F. Twilt, {\em Critical sets in parametric optimization}, Math. Program., 34 (1986), pp. 333-353. 
	
	\bibitem{JJT00} H.Th. Jongen, P. Jonker and F. Twilt, {\em Nonlinear Optimization in Finite Dimensions}, Kluwer Academic, Dordrecht, 2000.
	
	\bibitem{JS12} H.Th. Jongen and V. Shikhman, {\em Bilevel optimization: on the structure of the feasible set}, Math. Program.,  136 (2012), pp. 65-89.
	
	\bibitem{KBHP08} G. Kunapuli, K. Bennett, J. Hu and J-S. Pang, {\em Classification model selection via bilevel programming}, Optim. Meth. Softw., 23 (2008), pp. 475-489.
	
	
	\bibitem{LMYZZ20} R. Liu, P. Mu, X. Yuan, S. Zeng and J. Zhang, {\em A generic first-order algorithmic Framework for bi-level programming beyond lower-level singleton}, Proceedings of the 37th International Conference on Machine Learning, PMLR 119:6305-6315, 2020.
	
	\bibitem{LYZ08} G.S. Liu, J.J. Ye and J.P. Zhu, {\em Partial exact penaty for mathematical programs with equilibrium constraints}, Set-Valued Anal., 16 (2008), pp. 785-804.
	
	\bibitem{MPEC1} Z.-Q. Luo, J-S. Pang, and D. Ralph, {\em Mathematical Programs with Equilibrium Constraints}, Cambridge University Press, Cambridge, United Kingdom, 1996. 
	

	
	{\bibitem{MMZ20} P. Mehlitz, L.I. Minchenko, and A.B. Zemkoho, {\em A note on partial calmness for bilevel optimization problems with linearly structured lower level}, Optim. Lett. (2020), pp. 1-15. }
	
		\bibitem{M7599} J. Mirrlees, {\em The theory of moral hazard and unobservable behaviour-- part I}, Review of Economic Studies, 66 (1999), pp. 3-22.
	
	\bibitem{M18} B.S. Mordukhovich, {\em Variational Analysis and Applications},  Springer, Cham, 2018.
		
		\bibitem{M20} B.S. Mordukhovich, {\em Bilevel Optimization and Variational Analysis}, In: S. Dempe, A.B. Zemkoho (eds.)  Bilevel Optimization, Springer Optimization and its Applications, vol. 161, 2020.
	
	\bibitem{MNP12} B.S. Mordukhovich, N.M. Nam and H.M. Phan, {\em Variational analysis of marginal functions with applications to bilevel programming}, J. Optim. Theory Appl., 152 (2012), pp. 557-586.
	
	
	\bibitem{M51} T.S. Motzkin, {\em Two Consequences of the Transposition Theorem on Linear Inequalities}, Econometrica, 19 (1951), pp. 184-185.
	
	
	\bibitem{nie2017bilevel}
	J.~Nie, L.~Wang and J.J.~Ye,
	{\em  Bilevel polynomial programs and semidefinite relaxation methods},
	SIAM J.  Optim., 27 (2017), pp. 1728-1757.
	
\bibitem{nie2020bilevel}J.~Nie, L.~Wang, J.J.~Ye and S.~Zhong,
	{\em  A Lagrange multiplier expression method for bilevel polynomial programs}, SIAM J. Optim., 31 (2021), pp. 2368–2395.
	
	\bibitem{O90} J.V. Outrata, {\em On the numerical solution of a class of Stackelberg problems}, ZOR-Math. Methods Oper. Res., 34 (1990), pp. 255-277.
	
	
	
	\bibitem{MPEC2} J.V. Outrata, M. Kocvara and J. Zowe, {\em Nonsmooth Approach to Optimization Problems with Equlilibrium Constraints: Theory, Applications, and Numerical Results}, Kluwer Academic Publisher, Dordrect, The Netherlands, 1998.
	
	\bibitem{RW98} R.T. Rockafellar and R.J-B. Wets, {\em Variational Analysis}, Springer, Berlin, 1998.
	
	\bibitem{SS01} S. Scholtes and M. St\"{o}hr, {\em How stringent is the linear independence assumption for mathematical programs with complementarity constraints?}, Math. Oper. Res., 26 (2001), pp. 851–863.
	
	\bibitem{SIB97} K. Shimizu, Y. Ishizuka and J. Bard, {\em Nondifferentiable and Two-level Mathematical Programming}, Kluwer Academic, 1997.
	
	
	\bibitem{Ye04} J.J. Ye, {\em Nondifferentiable multiplier rules for optimization and bilevel optimization problems}, SIAM J. Optim., 15 (2004), pp. 252-274.
	
	\bibitem{Ye06} J.J. Ye, {\em Constraint qualifications and KKT conditions for bilevel programming problems}, Math. Oper. Res., 31 (2006), pp. 811-824.
	
	\bibitem{Ye11} J.J. Ye, {\em Necessary optimality  conditions for multiobjective bilevel programming problems}, Math. Oper. Res., 36 (2011), pp. 165-184.
	
	\bibitem{Ye20} J.J. Ye, {\em Constraint qualifications and optimality conditions in bilevel optimization},  In: S. Dempe, A.B. Zemkoho (eds.)  Bilevel Optimization, Springer Optimization and its Applications, vol. 161, 2020.
	
	\bibitem{YY97} J.J. Ye and X.Y. Ye, {\em  Necessary optimality conditions for optimization problems with variational inequality constraints}, Math. Oper. Res., 22 (1997), pp. 977-997.
	
	{\bibitem{YYZZ} J.J. Ye, X.M. Yuan, S. Zeng and J. Zhang, {\em Difference of convex algorithms for bilevel programs with applications in hyperparameter selection}, preprint, arXiv 2102.09006.}
	
		\bibitem{YZ18} J.J. Ye and J.C. Zhou, {\em  Verifiable sufficient conditions for the error bound property of second-order cone complementarity problems}, Math. Program., 171 (2018), pp. 361-395.
	
	\bibitem{YZ95} J.J. Ye and D.L. Zhu, {\em Optimality conditions for bilevel programming problems}, Optimization, 33 (1995), pp. 9-27.
	
	\bibitem{YZ10} J.J. Ye and D.L. Zhu, {\em New necessary optimality conditions for bilevel programs by combining the MPEC and value function approaches}, SIAM J. Optim., 20 (2010), pp. 1885-1905.
	
	
	
	\bibitem{YZZ97} J.J. Ye, D.L. Zhu and Q.J. Zhu, {\em Exact penalization and necessary optimality conditions for generalized bilevel programming problems}, SIAM J. Optim., 7 (1997), pp. 481-507.
	

	
	
\end{thebibliography}

\end{document}